\documentclass[11pt]{article}
\usepackage[T1]{fontenc}
\usepackage{lmodern}
\usepackage{fullpage}
\usepackage{mathrsfs}
\usepackage{authblk}
\usepackage{amssymb}
\usepackage{amsthm}
\usepackage{bbm}
\usepackage{setspace}
\usepackage{xspace}
\usepackage{multicol}
\usepackage{enumerate}
\usepackage{enumitem}
\usepackage{mathtools}
\usepackage{thm-restate}
\usepackage[hidelinks]{hyperref}
\usepackage{cleveref}

\emergencystretch=\maxdimen
\allowdisplaybreaks

\theoremstyle{plain}
\newtheorem{theorem}{Theorem}[section]
\newtheorem{lemma}[theorem]{Lemma}
\newtheorem{prop}[theorem]{Proposition}
\newtheorem{cor}[theorem]{Corollary}

\newtheorem*{claim*}{Claim}

\theoremstyle{definition}
\newtheorem{defn}[theorem]{Definition}
\newtheorem*{defn*}{Definition}

\newtheorem{remark}[theorem]{Remark}
\newtheorem{remarks}[theorem]{Remarks}
\newtheorem{question}[theorem]{Question}
\newtheorem*{question*}{Question}
\newtheorem{terminology}[theorem]{Terminology}

\newtheoremstyle{TheoremLub}
        {\topsep}{\topsep}                     
        {\itshape}                             
        {}                                     
        {\bfseries}                            
        {~\textnormal{(Lubarsky)}.}            
        { }                                    
        {\thmname{#1}\thmnote{ \bfseries #3}}  
    \theoremstyle{TheoremLub}
    \newtheorem{duplicateLub}{Lemma}

\newcommand{\tf}[1]{\textnormal{\textsc{#1}}\xspace}
\newcommand{\ff}[1]{\textnormal{\textrm{#1}}\xspace}

\newcommand{\Lhier}{\ensuremath{\mathbb{L}}\xspace}
\newcommand{\Lhieromega}{\tf{L}}

\DeclareMathOperator{\Defop}{\ff{Def}}
\DeclareMathOperator{\defop}{\ff{def}}
\DeclareMathOperator{\defclosure}{\mathfrak{D}}
\DeclareMathOperator{\defclosuresmall}{\mathfrak{E}}

\newcommand{\minusinf}{\textnormal{\textit{--Inf}}\xspace}
\newcommand{\pointwise}{\textnormal{``}}

\DeclareMathOperator{\setdiff}{\setminus}

\DeclareMathOperator{\downwards}{\downarrow}
\DeclareMathOperator{\mathand}{\mathrel{\wedge}}
\DeclareMathOperator{\mathor}{\mathrel{\vee}}

\DeclareMathOperator{\divline}{\hspace{0.05cm}\mid\hspace{0.05cm}}
\DeclareMathOperator{\bigdivline}{\hspace{0.05cm} \bigm\lvert \hspace{0.05cm}}
\DeclareMathOperator{\Bigdivline}{\hspace{0.05cm} \Bigm\lvert \hspace{0.05cm}}

\newcommand{\uminus}{\raisebox{.15\height}{\scalebox{0.75}{\ensuremath{-}}}}

\newcommand{\funop}[2]{\ensuremath \mathcal{#1}_{\mbox{\hspace{-4pt} \fontsize{5}{0}\ensuremath{#2}}}\xspace}

\DeclareMathOperator{\funpair}{\funop{F}{p}}
\DeclareMathOperator{\funint}{\funop{F}{\cap}}
\DeclareMathOperator{\fununion}{\funop{F}{\cup}}
\DeclareMathOperator{\fundiff}{\funop{F}{\setdiff}}
\DeclareMathOperator{\funtimes}{\funop{F}{\times}}
\DeclareMathOperator{\funimp}{\funop{F}{\rightarrow}}
\DeclareMathOperator{\funforall}{\funop{F}{\forall}}
\DeclareMathOperator{\fundom}{\funop{F}{d}}
\DeclareMathOperator{\funran}{\funop{F}{r}}
\DeclareMathOperator{\funabc}{\funop{F}{123}}
\DeclareMathOperator{\funacb}{\funop{F}{132}}
\DeclareMathOperator{\funeq}{\funop{F}{=}}
\DeclareMathOperator{\funin}{\funop{F}{\in}}

\newcommand{\rlz}{\Vdash_{\mathfrak{w}}}

\newcommand{\rlzt}{\Vdash_{\mathfrak{wt}}}

\newcommand{\rlztp}{\Vdash_{\mathfrak{wt}}^{\wp}}

\newcommand{\pair}{\textbf{p}}
\newcommand{\pairl}{\pair_0}
\newcommand{\pairr}{\pair_1}

\title{Constructing the Constructible Universe Constructively}
\date{}
\author[1]{Richard Matthews}

\affil[1]{Universit{\'{e}} Paris-Est Cr\'{e}teil, France}

\author[2]{Michael Rathjen}

\affil[2]{University of Leeds, United Kingdom}

\begin{document}

\maketitle

\begin{abstract}
We study the properties of the constructible universe, L, over intuitionistic theories. We give an extended set of fundamental operations which is sufficient to generate the universe over Intuitionistic Kripke-Platek set theory without Infinity. Following this, we investigate when L can fail to be an inner model in the traditional sense. Namely, we show that over Constructive Zermelo-Fraenkel (even with the Power Set axiom) one cannot prove that the Axiom of Exponentiation holds in L.
\end{abstract}

\section{Introduction}

In general, working with an arbitrary ``\emph{model of set theory}'' can be complicated. For example; the continuum could be arbitrarily large; there may be no definable global well-order or there may not be a clear, combinatorial structure to the universe. Instead, we may want to work with a ``\emph{nice}'' universe which has a clear fine structure that is well understood and permits a detailed analysis. Such a structure will usually also satisfy additional, useful axioms, for example the existence of a definable global well-order. The general framework for this is the concept of an \emph{inner model}. Over the standard set-theoretic framework of Zermelo-Fraenkel, \tf{ZF}, \tf{M} is said to be an inner model of $\tf{N} \models \tf{ZF}$ if:
\begin{itemize}
\item $\tf{M} \subseteq \tf{N}$,
\item $\tf{M} \models \tf{ZF}$,
\item $\tf{M} \cap \tf{Ord} = \tf{N} \cap \tf{Ord}$
\end{itemize}
and the simplest such example of this is G\"{o}del's \emph{Constructible Universe}, \tf{L}.

The Constructible Universe was developed by G\"{o}del in two influential papers, \cite{GodelLDefOp} and \cite{GodelLGodFun}, in the late 1930s in order to prove the consistency of the Axiom of Choice and the Generalised Continuum Hypothesis relative to \tf{ZF}. As with the standard $\tf{V}_\alpha$ hierarchy, one builds \tf{L} in stages where in place of the usual power set operation one instead uses a ``\emph{definable}'' version. The idea being that this thins out the universe one constructs by only accepting new sets which can be built from the previous sets in some definable way. Therefore, the main question one has to tackle is what is meant by ``definable''. That is, given $\tf{L}_\alpha$, what is $\tf{L}_{\alpha+1}$. 

There are three main ways to define this operation, all of which can be undertaken in the weak system of Kripke-Platek set theory: \begin{enumerate}
\item Syntactically; using the notion of a \emph{``definability operator''} so that $\tf{L}_{\alpha+ 1}$ is the collection of definable subsets of $\tf{L}_\alpha$. This is the original approach taken by G\"{o}del in \cite{GodelLDefOp} and was formalised in Kripke-Platek by Devlin in \cite{DevlinBook}.
\item By closure under what Barwise calls \emph{``fundamental operations''} or \emph{``G\"{o}del operations''}. This is the approach taken by G\"{o}del in \cite{GodelLGodFun} and further studied by Barwise in \cite{BarwiseAdmissibles} where he considered the theory of Kripke-Platek with urelements.
\item Using \emph{``rudimentary''} functions. This is a modified version of using fundamental operations which was developed\hspace{2pt} by Jensen and further explored by Mathias, leading to his weak system of Provi, the weakest known system in  which one can do both set forcing and build \tf{L}. The details of this theory can be found in \cite{MathiasBowlerRudimentary}.
\end{enumerate}

\noindent In this paper we address the question of inner models in intuitionistic set theories. As part of this, we first investigate a weak theory which is sufficient to construct \tf{L}. After this we question whether \tf{L} satisfies the conditions of being an inner model. \\

\noindent The intuitionistic approach to constructing \tf{L} was first undertaken by Lubarsky in \cite{LubarskyIntuitionisticL} under the assumption that \tf{V} satisfied Intuitionistic Zermelo-Fraenkel, \tf{IZF}. His approach was to show that the syntactic definition of the constructible universe still goes through in intuitionistic logic, with some minor modifications. The main obstacle one has to overcome is that the ordinals are no longer linearly ordered so one has to be more careful as to how one finds witnesses for the collection of definable subsets of some given set $X$. This also adds complications to proving the Axiom of Constructibility, that is proving that $\tf{V} = \tf{L}$ holds in \tf{L}. The main reason for this added difficulty is that, because it is unclear if \tf{L} should contain every ordinal under \tf{IZF}, on the face of it there is no reason to assume that 
\[
\bigcup_{\alpha \, \in \, \tf{Ord} \: \cap \: \tf{V}} \tf{L}_\alpha \quad = \quad \bigcup_{\alpha \, \in \, \tf{Ord} \: \cap \: \tf{L}} \tf{L}_\alpha.
\]
\noindent In order to circumvent this issue, Lubarsky proves the following lemma, which we shall reprove later in our weaker context of Intuitionistic Kripke-Platek set theory without infinity, $\tf{IKP}^{\minusinf}$.

\begin{duplicateLub}[\ref{thm: approx in L}]
For every ordinal $\alpha$ in \tf{V} there is an ordinal $\alpha^*$ in $\Lhieromega$ such that $\Lhieromega_\alpha = \Lhieromega_{\alpha^*}$.
\end{duplicateLub}

\noindent The syntactic approach has been further studied by Crosilla, and appears in the appendix to her PhD thesis, \cite{CrosillaPhD}. Here she shows that the construction can be carried out in a fragment of constructive set theory, which is equivalent to what we have defined as \tf{IKP}, by essentially the same proof as found in \cite{DevlinBook}. 

The third approach, via rudimentary functions, has also been explored in constructive contexts by Aczel \cite{AczelRudimentaryCST}. Here he defines the weak system of \emph{Rudimentary Constructive Set Theory} and shows that many of Jensen's techniques can be applied in this theory. \\

\noindent In this paper, we shall be interested in which axioms are sufficient to construct the constructible universe. Because the syntactic approach requires essential use of $\omega$ in order to work with arbitrarily long finite sequences, it is not the appropriate method to use in \tf{IKP} without infinity. Therefore, we shall adapt the second approach and use the fundamental operations. Adapting Barwise's method, we shall show that if one expands the collection of fundamental operations, then one can indeed construct \tf{L} over the weak system of \tf{IKP} without infinity. It should be noted that one could also consider urelements, as Barwise does. We have chosen not to undertake this study but this could be done without a significant amount of additional work. 

Following the structure of Barwise's original argument in Chapter II of \cite{BarwiseAdmissibles}, the main difference will be that we need to consider some additional logical operations. These are the operations of conjunction, implication and bounded universal quantification which are not treated in the classical case as distinct cases due to their equivalent definitions using disjunction, negation and bounded existential quantification. \\

\noindent Having constructed \tf{L}, we consider what happens when the universe is a model of Constructive Zermelo-Fraenkel set theory with Power Set, $\tf{CZF}_{\mathcal{P}}$. By combining a proof-theoretic ordinal analysis of Power Kripke-Platek with a notion of realizability from \cite{RathjenExistenceProperty} that implies truth, we shall see that Exponentiation need not hold in \tf{L}. Since Exponentiation follows from the axioms of Constructive Zermelo-Fraenkel set theory this means that \tf{CZF} cannot prove that \tf{L} is an inner model in the traditional sense.

\begin{restatable*}{theorem}{ExpFailsinLofCZF} \label{theorem:ExpFailsInLofCZFP}
$\tf{CZF}_\mathcal{P} \not\vdash \ff{Exp}^\tf{L}$. Thus, $\tf{CZF}_\mathcal{P}\not\vdash (\tf{CZF})^{\tf{L}}$.
\end{restatable*}

\noindent This result should be appreciated with the context of the analogous result over Power Kripke-Platek. In Proposition 6.47 of \cite{MathiasMacLane}, Mathias proves that Power Set, which is classically equivalent to Exponentiation, need not hold in \tf{L} if the background theory only satisfies Power Kripke-Platek. \\

\noindent So, where does this leave us in the search for intuitionistic inner models? We have that they need not satisfy every axiom of the base theory, it is far from clear if they contain every ordinal and it is even unclear if they contain a nice, combinatorial structure. Classically, there is an easy necessary and sufficient condition for a transitive class to be an inner model of \tf{ZF}. Moreover, this has the benefit of being expressible by a single sentence:

\begin{theorem}[\cite{JechSetTheory} 13.9]
A transitive class \tf{M} is an inner model of \tf{ZF} if and only if it is closed under the fundamental operations and is almost universal. That is, for every set $x$ which is a subset of $\tf{M}$, there is some $y \in \tf{M}$ for which $x \subseteq y$. 
\end{theorem}

\noindent We are unable to get such a strong connection in our intuitionistic context, due to the possible ill behaviour of our ordinals. On the other hand, we are able to get a slightly more modest result which still contains the essence of the original theorem. This is done by defining the notion of an \emph{external cumulative hierarchy} which is a notion encapsulating the idea that our model is iteratively constructed in the external universe via a class length sequence of set-sized approximations. For example,
\vspace{-5pt}
\[
\bigcup_{\alpha \in \tf{Ord} \cap \tf{V}} \tf{L}_\alpha.
\]

\begin{restatable*}{theorem}{ConditionforMtoModelIZF} \label{thm: condition for M to model IZF}
Suppose that \tf{V} is a model of \tf{IZF} and $\tf{M} \subseteq \tf{V}$ is a definable, transitive proper class with an external cumulative hierarchy. Then \tf{M} is a model of \tf{IZF} if and only if \tf{M} is closed under the fundamental operations and is almost universal.
\end{restatable*}

\section{Intuitionism}

When working intuitionistically, it is important to specify the axioms correctly in order to avoid inadvertently being able to deduce the Law of Excluded Middle. For example, we have the following classical results which can be shown to go through over very weak theories.

\begin{theorem}[\cite{AczelRathjenCST} 10.4.1]
The Foundation Scheme implies the Law of Excluded Middle. 
\end{theorem}

\begin{theorem}[Diaconescu, \cite{DiaconescuConstructiveAC}] \index{Theorem! Diaconescu}
The Axiom of Choice implies the Law of Excluded Middle.
\end{theorem}

\noindent Avoiding these principles leads to the following axiomatisations of \emph{Intuitionistic Zermelo-Fraenkel}, \tf{IZF}, and \emph{Intuitionistic Kripke-Platek}, \tf{IKP}. It is proven in \cite{FriedmanScedrovDefinableWitnesses} that, intuitionistically, the Replacement Scheme does not imply the Collection Scheme, even if one assumes Dependent Choice holds along with every classically true $\Sigma_1$-sentence in the language of \tf{ZF}. Moreover, Replacement alone is insufficient to do many of the constructions we need, so we shall formulate \tf{IZF} with the Collection Scheme. Secondly, in theories with restricted Separation it is beneficial to strengthen the Axiom of Infinity to asserting that there is an inductive set which is contained in every other inductive set, leading to the Axiom of Strong Infinity. Even over \tf{IKP} this is known to be equivalent to the Axiom of Infinity (for example see Proposition 4.7 of \cite{AczelRathjenTechnicalReport} where this equivalence is proven over \tf{CZF}) however this is non-trivial, so we include the stronger version here for simplicity.

\begin{defn}[Axiom of Strong Infinity]
$ \exists a \: (\ff{Ind}(a) \mathand \forall b \: (\ff{Ind}(b) \rightarrow \forall x \in a (x \in b)))$ where $\ff{Ind}(a)$ is an abbreviation for $\emptyset \in a \mathand \forall x \in a (x \cup \{ x \} ) \in a $.
\end{defn}

\begin{defn}
Let \tf{IZF} denote the theory whose underlying logic is intuitionistic and whose axioms are; Extensionality, Empty Set, Pairing, Unions, Power Set, Strong Infinity, Set Induction Scheme, Separation Scheme and Collection Scheme.
\end{defn}

\begin{defn}
Let \tf{IKP} denote the theory whose underlying logic is intuitionistic and whose axioms are; Extensionality, Empty Set, Pairing, Unions, Strong Infinity, Set Induction Scheme, $\Sigma_0$-Separation Scheme and $\Sigma_0$-Collection Scheme.
\end{defn}

\noindent Because we are interested in which axioms are necessary to construct the constructible universe, in \Cref{sec: fundamental operations} we shall take care to differentiate between \tf{IKP} \emph{without infinity}, which we will call $\tf{IKP}^{\minusinf}$, and \tf{IKP}.\footnote{It is worth noting that some authors, for instance Barwise in \cite{BarwiseAdmissibles}, do not include the Axiom of Infinity in their formulation of \tf{KP} but instead use either \tf{KPI} or $\tf{KP}\omega$ to denote the theory with Infinity. Moreover, other authors, such as Mathias in \cite{MathiasWeakSystems} formulate \tf{KP} with only $\Pi_1$-Foundation. All the results in this paper would go through with this weakening but we choose to include full Foundation for simplicity.}  We shall therefore regularly refer to \tf{IKP} as $\tf{IKP}^{\minusinf} + \emph{Strong Infinity}$ just to make it clear when Strong Infinity is being assumed. 

\smallskip

\noindent While Diaconescu's theorem tells us that the Axiom of Choice\footnote{in the form that for every set $A$ and function $F$ with domain $A$ satisfying $\forall x \in A \: \exists y \in F(x)$ there is a function $f$ with domain $A$ such that $\forall x \in A \: f(x) \in F(x)$} results in instances of Excluded Middle, this is not the case for all Choice Principles. For example, one can add Dependent Choice (for example to Bishop's constructive mathematics \cite{BishopBridgesConstructiveAnalysis})  and, by a result from Chapter 8 of \cite{BellSetTheory}, a version of Zorn's Lemma holds in every Heyting-valued model over a model of \tf{ZFC}. As a second remark, one can see that many of the basic properties of $\tf{KP}^{\minusinf}$ can still be deduced in $\tf{IKP}^{\minusinf}$. We note a few such theorems here, the proofs of which are given in Chapter 19 of \cite{AczelRathjenCST}.

\begin{prop} In $\tf{IKP}^{\minusinf}$ one can prove that:
\begin{itemize} \setlength \itemsep{0pt}
\item Let $\langle a, b \rangle \coloneqq \{ \{a\}, \{a, b\} \}$. Then $\forall a, b, c, d \: (\langle a, b \rangle = \langle c, d \rangle \rightarrow (a = c \mathand b = d ))$.
\item $\forall a, b \: \exists c \: (c = a \times b)$.
\end{itemize}
\end{prop}

\noindent Another important fact about $\tf{IKP}^{\minusinf}$ is that we can still deduce the Reflection and Strong Collection Principles for a larger class of formulae, namely the class of $\Sigma$-formulae where

\begin{defn}
The class of $\Sigma$-formulae is the smallest class containing the atomic formulae which is closed under $\wedge$, $\vee$, $\rightarrow$, the bounded quantifiers $\forall x \in a$ and $\exists x \in a$ and the unbounded existential quantifier $\exists x$.
\end{defn}

\begin{theorem}[$\Sigma$-Reflection Principle]
For every $\Sigma$-formula $\varphi(v)$ we have 
\[
\tf{IKP}^{\minusinf} \vdash \varphi(v) \leftrightarrow \exists a \: \varphi^{(a)}(v),
\]
where $\varphi^{(a)}(v)$ is the result of replacing each unbounded quantifier in $\varphi$ with bounded quantification over $a$.
\end{theorem}

\begin{theorem}[$\Sigma$-Strong Collection Principle]
For every $\Sigma$-formula $\varphi(x, y, z)$ the following is a theorem of $\tf{IKP}^{\minusinf}$: For every set $z$, if $\forall x \in a \: \exists y \: \varphi(x, y, z)$ then there is a set $b$ such that 
\[
\forall x \in a \: \exists y \in b \: \varphi(x, y, z) \mathand \forall y \in b \: \exists x \in a \: \varphi(x, y, z).
\]
\end{theorem}

\noindent One final, particularly powerful tool that is provable in $\tf{IKP}^{\minusinf}$ is definitions by $\Sigma$-recursion. This allows one to prove that $\Sigma$-definable operations, such as the transitive closure of a set, are provably total.

\begin{theorem}[$\Sigma$-Recursion]
If $G$ is a total $(n+2)$-ary $\Sigma$-definable class function, then there exists a total $(n+1)$-ary $\Sigma$-definable class function $F$ such that 
\[
\tf{IKP}^{\minusinf} \vdash \forall \overrightarrow{x} \forall y \big( F(\overrightarrow{x}, y) = G(\overrightarrow{x}, y, \{ \langle z, F(\overrightarrow{x}, z) \rangle \divline z \in y \} ) \big).
\]
\end{theorem}

\noindent However we are not able to deduce the $\Delta$-Separation Principle. This is because the proof makes essential use of the fact that classically $ \forall x \in a (\varphi \rightarrow \psi)$ holds if and only if $\forall x \in a (\psi(x) \mathor \neg \varphi(x))$ holds, which is not intuitionistically  provable. \\

\noindent An important concept is the notion of an ordinal. This is traditionally defined as a transitive set which is well-ordered by the $\in$ relation, however this is not appropriate in our setting. Therefore, we shall follow the definition given in I.3.6 of Barwise \cite{BarwiseAdmissibles}, which is equivalent to the standard one over \tf{KP}:

\begin{defn} \index{Ordinal} 
An \emph{ordinal} is a transitive set of transitive sets.
\end{defn}

\noindent Let $0$ denote the ordinal $\emptyset$, which is the unique set given by the Axiom of Empty Set, and $1$ the ordinal $\{0\}$. An important class of ordinals is the class of \emph{Truth Values} 
\[
\Omega \coloneqq \{ x \divline x \subseteq 1 \}.
\]
\noindent $\Omega = \{0, 1 \}$ is the assertion that every statement is either true or false, which entails the Law of Excluded Middle. Therefore, in general we may have other truth values even if they can't necessarily be determined. For example, these can be ordinals of the form 
\[
\{ 0 \in 1 \divline \varphi \}
\]
\noindent where $\varphi$ is a formula for which we can neither determine $\varphi$ nor $\neg \varphi$. The complexity introduced by having more than 2 truth values means that the natural ordering assigned to the ordinals becomes tree-like rather than a simple linear order, as can be seen by the following proposition.

\begin{prop} \,
\[
\tf{IKP}^{\minusinf} \vdash \forall \alpha \in \tf{Ord} (0 \in \alpha + 1) \longrightarrow \textit{Law of Excluded Middle for } \Sigma_0 \textit{-formulae}.
\]
\end{prop}

\begin{proof}
Let $\varphi$ be a $\Sigma_0$-formula and consider the ordinal truth value $\alpha = \{ 0 \in 1 \divline \varphi\}$. If $0$ were an element of $\alpha + 1 = \alpha \cup \{\alpha\}$ this would imply that either $0 \in \alpha$ or $0 = \alpha$. The first of these scenarios implies that $\varphi$ holds while the second one gives us $\neg \varphi$, and thus $\varphi \mathor \neg \varphi$ must be true.
\end{proof}

\section{Intermediate Theories}

\noindent In this section we include several more intuitionistic theories intermediate between \tf{IKP} and \tf{IZF}. The first one extends \tf{IKP} by adding power sets.

In order to do this, we begin by defining the class of $\Sigma_0^{\mathcal{P}}$-formulae. This extends the class $\Sigma_0$ by adding in bounded universal and existential subset quantifiers. Let $\exists x \subseteq y$ and $\forall x \subseteq y$ be abbreviations for $\exists x (x \subseteq y \mathand \dots)$ and $\forall x (x \subseteq y \rightarrow \dots)$ respectively. 

\begin{defn}
The class of $\Sigma^{\mathcal{P}}_0$-formulae is the smallest class which contains the atomic formulae and is closed under $\mathand, \mathor, \rightarrow$ and the quantifiers
\[
\forall x \in a, \: \exists x \in a, \: \forall x \subseteq a, \: \exists x \subseteq a.
\] 
\end{defn}

\begin{defn}
Let $\tf{IKP}(\mathcal{P})$ denote the theory whose underlying logic is intuitionistic and whose axioms are; Extensionality, Empty Set, Pairing, Unions, Power Set, Strong Infinity, Set Induction Scheme, $\Sigma^{\mathcal{P}}_0$-Separation Scheme and $\Sigma^{\mathcal{P}}_0$-Collection Scheme.
\end{defn}

\begin{remark}[Mathias]
We observe here that $\tf{KP}(\mathcal{P})$ is different from \tf{KP} plus Power Set. For example, $\tf{L}_{\aleph^{\tf{L}}_\omega} \models \tf{KP} + \textnormal{ Power Set}$ but does not model $\tf{KP}(\mathcal{P})$.
\end{remark}

\noindent The next theories we shall discuss are variants of Aczel's \emph{Constructive Zermelo-Fraenkel}, an in-depth study of this theory can be found in \cite{AczelRathjenCST}. For these we need the axioms of \emph{Strong Collection} and \emph{Subset Collection}:

\begin{defn} \,

\noindent (Strong Collection). For any formula $\varphi$ and set $a$,
\[ 
\forall x \in a \: \exists y \: \varphi(x,y) \longrightarrow \exists b \Big( \forall x \in a \: \exists y \in b \: \varphi(x, y) \mathand \forall y \in b \: \exists x \in a \: \varphi(x, y) \Big).
\]
(Subset Collection). For any formula $\varphi$ and sets $a, b$,
\[
\exists c \forall u \bigg( \forall x \in a \: \exists y \in b \: \varphi(x, y, u) \rightarrow \exists d \in c \: \Big( \forall x \in a \: \exists y \in d \: \varphi(x, y, u) \mathand \forall y \in d \: \exists x \in a \: \varphi(x, y, u) \Big) \bigg).
\]
\end{defn}

\begin{defn} \,
\begin{itemize}
\item Let $\tf{CZF}^{\uminus}$ denote \tf{IKP} plus Strong Collection,
\item Let $\tf{CZF}$ denote $\tf{CZF}^{\uminus}$ plus Subset Collection,
\item Let $\tf{CZF}_\mathcal{P}$ denote $\tf{CZF}^{\uminus}$ plus Power Set.
\end{itemize}
\end{defn}

\medskip

\begin{remarks} \,
\begin{itemize} 
\item Clearly $\tf{CZF}$ is a subtheory of $\tf{CZF}_\mathcal{P}$ and the former is known to prove the Axiom of Exponentiation (see Theorem 5.1.2 of \cite{AczelRathjenCST}). Therefore, in this context, Subset Collection can be seen as saying that a weak form of Power Set holds.
\item One can easily show, by induction on the complexity of $\Sigma^{\mathcal{P}}_0$-formulae, that $\tf{CZF}_\mathcal{P}$ proves every instance of $\Sigma^{\mathcal{P}}_0$-Separation. Therefore, $\tf{IKP}(\mathcal{P})$ is a subtheory of $\tf{CZF}_\mathcal{P}$.
\item An alternative way to view $\tf{CZF}_{\mathcal{P}}$ is as \tf{IZF} with Full Separation weakened to only $\Sigma_0^{\mathcal{P}}$-Separation.
\item The theories $\tf{CZF}^{\uminus}$ and \tf{CZF} are known to be equiconsistent with \tf{KP}. However, $\tf{CZF}$ plus Full Excluded Middle results in \tf{ZF}. On the other hand, as proven in Theorem 15.1 of \cite{RathjenCalculusofConstructions}, $\tf{CZF}_\mathcal{P}$ is much stronger than \tf{KP} and even stronger than Zermelo set theory.
\end{itemize}
\end{remarks}

\section{The Fundamental Operations} \label{sec: fundamental operations}

We now proceed to construct the constructible universe, \tf{L}. In order to do this, given $\tf{L}_\alpha$ we would like $\tf{L}_{\alpha+1}$ to be the collection of ``\emph{definable subsets}'' of $\tf{L}_\alpha$ for some appropriate definition of definable. Following Barwise, we will choose to do this by defining a series of $\Sigma$-operations, each of which will comprise of two arguments. These operations will then be used to generate the constructible sets. The main benefit to this approach over the syntactic one is that it will avoid assuming any form of the Axiom of Infinity.

The operations we will take are the same as in \cite{BarwiseAdmissibles} except for the addition of $\funint$, $\funimp$ and $\funforall$. This is to allow us to prove \Cref{lemma:BarwisetFormula} that the Axiom of Separation for a $\Sigma_0$-formula can be obtained by a sequence of such operations. This will be achieved by an induction on the construction of the class of $\Sigma_0$ formulae and in order to generalise this to the intuitionistic setting we will need to consider the three additional cases of conjunction, implication and bounded universals that are not treated by Barwise. Conjunctions will be relatively straight forward and we observe here that, over classical logic, intersections are deducible from the original operations via the identity 
\[
x \cap y = x \setdiff (x \setdiff y),
\]
\noindent but this equivalence does not hold intuitionistically. Since implication cannot be deduced from the other logical connectives we need an operation that specifically does this. To be more specific, suppose we have sequences of fundamental operations $\mathcal{F}_\varphi$ and $\mathcal{F}_\psi$ such that $\mathcal{F}_\varphi(a) = \{ z \in a \divline \varphi(z) \}$ and $\mathcal{F}_\psi(a) = \{ z \in a \divline \psi(z)\}$ and we wanted a sequence $\mathcal{G}$ such that $\mathcal{G}(a) = \{ z \in a \divline \varphi(z) \rightarrow \psi(z) \}$. This set is equal to $\{ z \in a \divline z \in \mathcal{F}_\varphi(a) \rightarrow z \in \mathcal{F}_\psi(a) \}$ and therefore we need an operation that, given sets $a$, $b$ and $c$, returns $\{ z \in c \divline z \in a \rightarrow z \in b \}$ (where $c$ is needed to ensure this operation results in a set) which is what $\funimp$ does.

The final operation is to allow us to find fundamental sequences giving bounded universal quantification. This will be done by adding the operation $\funforall$ which essentially says that, given a set relation $x$ the collection of preimages, $\{ \{ u \divline \langle z, u \rangle \in x \} \divline z \in \ff{dom}(x) \}$ forms a set. We remark here that the same operation is included in the axioms of Provi by Mathias, \cite{MathiasBowlerRudimentary}. We begin by first fixing some notation:

\begin{defn} 
For $x$ an ordered pair, $y$ a set of ordered pairs and $z$ a set, define
\begin{itemize} \setlength \itemsep{1pt}
\item $1^{st}(x) = a$ iff $\exists u \in x \: \exists b \in u \: (x = \langle a, b \rangle)$, 
\item $2^{nd}(x) = b$ iff $\exists u \in x \: \exists a \in u \: (x = \langle a, b \rangle)$,
\item $y \pointwise \{z\} \coloneqq \{ u \divline \langle z, u \rangle \in y\}.$
\end{itemize}
\end{defn}

\begin{defn}\label{operations} \index{Fundamental Operations}
The \emph{fundamental operations} are as follows:
\begin{itemize}[labelindent=\parindent, leftmargin=*, widest=$(\funabc)$, align=left] \setlength \itemsep{1pt}
\item[$(\funpair)$] $\funpair(x, y) \coloneqq \{ x, y\},$
\item[$(\funint)$] $\funint(x, y) \coloneqq x \cap \bigcap y,$
\item[$(\fununion)$] $\fununion(x, y) \coloneqq \bigcup x,$
\item[$(\fundiff)$] $\fundiff(x, y) \coloneqq x \setdiff y,$
\item[$(\funtimes)$] $\funtimes(x, y) \coloneqq x \times y,$
\item[$(\funimp)$] $\funimp(x, y) \coloneqq x \cap \{ z \divline \textit{y is an ordered pair} \, \mathand \, (z \in 1^{st}(y) \rightarrow z \in 2^{nd}(y)) \},$
\item[$(\funforall)$] $\funforall(x, y) \coloneqq \{ x \pointwise \{z\} \divline z \in y \} = \{ \{ u \divline \langle z, u \rangle \in x \} \divline z \in y \} ,$
\item[$(\fundom)$] $\fundom(x, y) \coloneqq \ff{dom}(x) = \{ 1^{st}(z) \divline z \in x \mathand \textit{z is an ordered pair}\},$
\item[$(\funran)$] $\funran(x, y) \coloneqq \ff{ran}(x) = \{ 2^{nd}(z) \divline z \in x \mathand \textit{z is an ordered pair}\},$
\item[$(\funabc)$] $\funabc(x, y) \coloneqq \{ \langle u, v, w \rangle \divline \langle u, v \rangle \in x \mathand w \in y \},$
\item[$(\funacb)$] $\funacb(x, y) \coloneqq \{ \langle u, w, v \rangle \divline \langle u, v \rangle \in x \mathand w \in y \},$
\item[$(\funeq)$] $\funeq(x, y) \coloneqq \{ \langle v, u \rangle \in y \times x \divline u = v \},$
\item[$(\funin)$] $\funin(x, y) \coloneqq \{ \langle v, u \rangle \in y \times x \divline u \in v \}.$
\end{itemize}
\end{defn}

\begin{remark} In order to simplify later notation, we shall let $\mathcal{I}$ be the obvious finite set indexing the above operations.
\end{remark}

\noindent Note that we form n-tuples inductively as 
\[
\langle x_3, x_2, x_1 \rangle \coloneqq \langle x_3, \langle x_2, x_1 \rangle \rangle,
\]
\noindent and therefore $\ff{ran}( \{ \langle x_3, x_2, x_1 \rangle \} ) = \{ \langle x_2, x_1 \rangle \}$.  
The next lemma is adapted from Lemma II.6.1 of \cite{BarwiseAdmissibles} which will give us that any instance of Bounded Separation can be written as a sequence of fundamental operations. Note that in the lemma the order of the variables has been inverted because this turns out to simplify the arguments.

\begin{lemma} \label{lemma:BarwisetFormula} For every $\Sigma_0$-formula $\varphi(x_1, \dots, x_n)$ with free variables among $x_1, \dots, x_n$, there is a term $\funop{F}{\varphi}$ built up from the operations in \Cref{operations} such that 
\[
\tf{IKP}^{\minusinf} \vdash \funop{F}{\varphi}(a_1, \dots, a_n) = \{ \langle x_n, \dots, x_1 \rangle \in a_n \times \ldots \times a_1 \divline \varphi(x_1, \dots, x_n) \}.
\]
\end{lemma}

\begin{proof}
As in Barwise, we will call a formula $\varphi(x_1, \dots, x_n)$ a \emph{termed-formula}, or \emph{t-formula}, if there is a term $\funop{F}{\varphi}$ built from the fundamental operations such that the conclusion of the lemma holds. We shall then proceed by induction on $\Sigma_0$-formulae to show that every such formula is a t-formula. Using the proof of \cite{BarwiseAdmissibles}, Lemma II.6.1, it only remains to consider the following cases: \begin{enumerate}[label=(\roman*.)]
\item \label{thm: P and Q is t-form} If $\varphi(x_1, \dots, x_n)$ and $\psi(x_1, \dots, x_n)$ are t-formulae then so is $\varphi(x_1, \dots, x_n) \mathand \psi(x_1, \dots, x_n)$.
\item \label{thm: P implies Q is t-form} If $\varphi(x_1, \dots, x_n)$ and $\psi(x_1, \dots, x_n)$ are t-formulae then so is $\varphi(x_1, \dots, x_n) \rightarrow \psi(x_1, \dots, x_n)$.
\item \label{thm: for all x in b is t-form} If $\psi(x_1, \dots, x_{n + 1})$ is a t-formula and $\varphi(x_1, \dots, x_n, b)$ is $\forall x_{n + 1} \in b \: \psi(x_1, \dots, x_{n + 1})$, where $b$ is an arbitrary set that does not appear in $\{x_1, \dots, x_n \}$, then $\varphi$ is a t-formula.
\item \label{thm: for all xi in xj is t-form} If $\psi(x_1, \dots, x_{n + 1})$ is a t-formula and $\varphi(x_1, \dots, x_n)$ is $\forall x_{n + 1} \in x_j \: \psi(x_1, \dots, x_{n + 1})$, then $\varphi$ is a t-formula.
\end{enumerate}

\noindent \textbf{Case} \ref{thm: P and Q is t-form}: Let $\funop{F}{\varphi}$ and $\funop{F}{\psi}$ witness that $\varphi$ and $\psi$ are t-formulae. First note that for a set $z$, $\funpair(z, z) = \{ z \}$. Then, 
\begin{align*}
\funint \bigg( \funop{F}{\varphi}(a_1, \dots, a_n), \funpair  & \Big( \funop{F}{\psi}(a_1, \dots, a_n), \funop{F}{\psi}(a_1, \dots, a_n) \Big) \bigg) \\
& = \funop{F}{\varphi}(a_1, \dots, a_n) \hspace{0.05cm} \cap \hspace{0.05cm} {\displaystyle \bigcap} \big\{ \funop{F}{\psi}(a_1, \dots, a_n) \big\} \\
& = \funop{F}{\varphi}(a_1, \dots, a_n) \cap \funop{F}{\psi}(a_1, \dots, a_n).
\end{align*}
\noindent Therefore, we can define $\funop{F}{\varphi \mathand \psi}(a_1, \dots, a_n)$ as 
\[
\funint \bigg( \funop{F}{\varphi}(a_1, \dots, a_n), \funpair \Big( \funop{F}{\psi}(a_1, \dots, a_n), \funop{F}{\psi}(a_1, \dots, a_n) \Big) \bigg).
\]

\noindent \textbf{Case} \ref{thm: P implies Q is t-form}: Let $\funop{F}{\varphi}$ and $\funop{F}{\psi}$ witness that $\varphi$ and $\psi$ are t-formulae. Then 
\begin{align*}
\{ \langle x_n, \dots, x_1 \rangle \in a_n \times & \ldots \times a_1 \divline \varphi(x_1, \dots, x_n) \rightarrow \psi(x_1, \dots, x_n) \} \\ 
& = (a_n \times \dots \times a_1) \cap \{ z \divline z \in \funop{F}{\varphi}(a_1, \dots, a_n) \rightarrow z \in \funop{F}{\psi}(a_1, \dots, a_n)\}.
\end{align*}
\noindent Also, 
\[
\langle x, y \rangle = \{ \{ x \}, \{ x, y \} \} = \funpair \big( \funpair(x, x), \funpair(x, y) \big)
\]
\noindent and $a_n \times \ldots \times a_1$ can be defined by repeated used of $\funtimes$ so we can use these constructions. Thus, the above can be expressed as 
\[
\funimp \Big( a_n \times \ldots \times a_1, \: \Big\langle \funop{F}{\varphi}(a_1, \dots, a_n), \funop{F}{\psi}(a_1, \dots, a_n) \Big\rangle \Big),
\]
\noindent giving the required construction of $\mathcal{F}_{\varphi \rightarrow \psi}.$ 

\smallskip

\noindent \textbf{Case} \ref{thm: for all x in b is t-form}: Let $\varphi(x_1, \dots, x_n, b) \equiv \forall x_{n + 1} \in b \: \psi(x_1, \dots, x_{n + 1})$ and let $\funop{F}{\psi}$ witness that $\psi$ is a t-formula. Then
\begin{eqnarray*}
\funforall \Big( \funop{F}{\psi}(a_1, \dots, a_n, b), b \Big) \hspace{-2pt} & = & \Big\{ \funop{F}{\psi}(a_1, \dots, a_n, b) \pointwise \{z \} \bigdivline z \in b \Big\} \\
& = & \bigg \{ \Big \{ w \bigdivline \langle z, w \rangle \in \funop{F}{\psi}(a_1, \dots, a_n, b) \Big \} \Bigdivline z \in b \bigg \} \\
& = & \bigg \{ \Big \{ \langle x_n, \dots x_1 \rangle \bigdivline \langle z, x_n, \dots, x_1 \rangle \in \funop{F}{\psi}(a_1, \dots, a_n, b ) \Big \} \Bigdivline z \in b \bigg \} \\
& = & \Big \{ \ff{ran}( \funop{F}{\psi}(a_1, \dots, a_n, \{z \})) \bigdivline z \in b \big \}.
\end{eqnarray*} 

\noindent Therefore $\funop{F}{\varphi}(a_1, \dots, a_n, b)$ can be expressed as 
\begin{eqnarray*}
\Big \{ \langle x_n, \dots, x_1 \rangle \hspace{-0.25cm} & \in & \hspace{-0.25cm} a_n \times \ldots \times a_1 \divline \forall x_{n + 1} \in b \: \psi(x_1, \dots x_{n+1}) \Big \} \\
& = & \hspace{-0.25cm} (a_n \times \ldots \times a_1) \cap \Big \{ w \bigdivline \forall x_{n + 1} \in b \; \langle x_{n + 1}, w \rangle \in \funop{F}{\psi} \big ( a_1, \dots, a_n, \{x_{n + 1}\} \big ) \Big \} \\
& = & \hspace{-0.25cm} (a_n \times \ldots \times a_1) \cap \bigcap  \big \{ \ff{ran}(\funop{F}{\psi}(a_1, \dots, a_n, \{x_{n + 1} \})) \divline x_{n + 1} \in b \big \} \\
& = & \hspace{-0.25cm} \funint \Big ( a_n \times \ldots \times a_1, \: \funforall \big ( \funop{F}{\psi} (a_1, \dots, a_n, b), b \big ) \Big ).
\end{eqnarray*} 

\noindent \textbf{Case} \ref{thm: for all xi in xj is t-form}: Let $\varphi(x_1, \dots, x_n) \equiv \forall x_{n + 1} \in x_j \: \psi(x_1, \dots, x_{n + 1})$. Then 
\[
\{ \langle x_n, \dots, x_1 \rangle \in a_n \times \ldots \times a_1 \divline \forall x_{n + 1} \in x_j \: \psi(x_1, \dots, x_{n + 1}) \} 
\]
\noindent is equal to the following set;
\[
\{ \langle x_n, \dots, x_1 \rangle \in a_n \times \ldots \times a_1 \divline \forall x_{n + 1} \in \bigcup a_j \: ( x_{n + 1} \in x_j \rightarrow \psi(x_1, \dots, x_{n + 1}) ) \}. 
\]
\noindent So, taking $\vartheta(x_1, \dots, x_n, \bigcup a_j) \equiv \forall x_{n + 1} \in \bigcup a_j \: ( x_{n + 1} \in x_j \rightarrow \psi(x_1, \dots, x_{n + 1}) ) $, $\varphi$ is a t-formula by cases \ref{thm: P implies Q is t-form} and \ref{thm: for all x in b is t-form} and the fact that if two formulae are provably equivalent in $\tf{IKP}$ and one is a t-formula then so is the other.\footnote{This is statement (b) in the proof of Lemma II.6.1 of \cite{BarwiseAdmissibles}}
\end{proof}

\begin{theorem}\label{thm: bounded sep}
For any $\Sigma_0$-formula $\varphi(x_1, \dots, x_n)$ with free variables among $x_1, \dots x_n$, there is a term $\funop{G}{\varphi}$ of $n$ arguments built from the operations defined in \Cref{operations} such that: 
\[
\tf{IKP}^{\minusinf} \vdash \funop{G}{\varphi}(a, x_1, \dots, x_{i - 1}, x_{i + 1}, \dots, x_n) = \{ x_i \in a \divline \varphi(x_1, \dots, x_n) \}.
\]
\end{theorem}

\begin{proof}
This follows easily from our lemma since if $\funop{F}{\varphi}$ is the term built in the previous lemma such that 
\[
\tf{IKP}^{\minusinf} \vdash \funop{F}{\varphi}(a_1, \dots, a_n) = \{ \langle x_n, \dots, x_1 \rangle \in a_n \times \ldots \times a_1 \divline \varphi(x_1, \dots, x_n) \}
\]
\noindent then our required set can be built from $\funop{F}{\varphi}(\{x_1\}, \dots, \{x_{i - 1} \}, a_i, \{ x_{i + 1} \}, \dots \{ x_n \})$ by using $\funran$ $n - i$ times and then $\fundom$. 
\end{proof}

\section{Defining Definability}

Having defined the fundamental operations, in this section we shall give a \emph{definability operator}. The idea will be that the definable subsets of $b$ are those sets which can be constructed from $b$ using the fundamental operations. We shall then discuss some of the basic properties one can deduce from this definition and show that the model which one constructs satisfies $\tf{IKP}^{\minusinf}$. To conclude, we will end this section by mentioning other definability operators.

\begin{defn} 
For a set $b$:
\begin{itemize}
\item $\defclosuresmall(b) \coloneqq b \cup \{ \funop{F}{i}(x, y) \divline x, y \in b \mathand i \in \mathcal{I} \}$,
\item $\defclosure(b) \coloneqq \defclosuresmall(b \cup \{b\}).$
\end{itemize}
\end{defn}

\noindent The following proposition is then provable using $\Sigma$-Collection and the Axiom of Unions in $\tf{IKP}^{\minusinf}$.

\begin{prop}
$\tf{IKP}^{\minusinf} \vdash \forall b \: \exists x \: (x = \defclosuresmall(b)).$
\end{prop}

\noindent Observe that, with this construction, $\defclosure(b)$ need not be transitive even if $b$ is a transitive set. However, following Theorem II.6.4 of \cite{BarwiseAdmissibles}, one could easily extend the list of fundamental functions to do this. Namely, one could add the following functions, which we state with three variables for simplicity:
\begin{itemize}
\begin{multicols}{2}
\item $\funop{G}{0}(x, y, z) \coloneqq \langle x, y \rangle$, 
\item $\funop{G}{2}(x, y, z) \coloneqq \langle x, y, z \rangle$,
\item $\funop{G}{4}(x, y, z) \coloneqq \{ x, \langle y, z \rangle \}$,

\item $\funop{G}{1}(x, y, z) \coloneqq \{ u \divline \langle y, u \rangle \in x \}$,
\item $\funop{G}{3}(x, y, z) \coloneqq \langle x, z, y \rangle$,
\item $\funop{G}{5}(x, y, z) \coloneqq \{ x, \langle z, y \rangle \}$.
\end{multicols}
\end{itemize}

\begin{prop}
For any set $b$, if $b$ is transitive then so is 
\[
b \cup \{ \funop{F}{i}(x. y) \divline x, y \in b \mathand i \in \mathcal{I} \} \cup \{ \funop{G}{i}(x, y, z) \divline x, y, z \in b \mathand i \in 6 \}.
\]
\end{prop}

\noindent However, apart from making the definable closure into a transitive set, there is no obvious benefit to using this expanded list. Moreover, each of the $\funop{G}{i}$ can be easily built from at most five applications of the original operations, for example, 
\[
\{ u \divline \langle y, u \rangle \in x \} = \funran \bigg( \funint \bigg( x \: , \: \funpair \Big( \funtimes \big( \funpair(y,y), \funran(x) \big), \funtimes \big( \funpair(y,y), \funran(x) \big) \Big) \bigg) \bigg).
\]
Therefore, we choose to just use the original list as these are the crucial operations. \\

\noindent We would now want to define another operation $\Defop(b)$ to be the closure of $b$ under our fundamental operations, that is 
\[
\Defop(b) \coloneqq \bigcup_{n \in \omega} \defclosure^n(b).
\]
\noindent This would have the added benefit that if we defined $\Lhieromega_\alpha \coloneqq {\displaystyle\bigcup}_{\beta \in \alpha} \Defop(\Lhieromega_\beta)$, then for each ordinal $\alpha$, $\Lhieromega_\alpha$ would be a transitive model of Bounded Separation. However, this definition requires the Axiom of Infinity, which we are not initially assuming. Therefore, to begin with, we will just use $\defclosure$ to define our universe and use the different script $\Lhier$ to differentiate between the two notions. The relationship between $\Lhier$ and $\Lhieromega$ in the presence of Strong Infinity will be discussed in \Cref{thm: with inf vs without}.

\begin{defn} For $\alpha$ an ordinal, ${ \displaystyle \Lhier_\alpha \coloneqq \bigcup_{\beta \in \alpha} \defclosure(\Lhier_\beta)}$ and ${ \displaystyle \Lhier \coloneqq \bigcup_{\alpha \in \tf{Ord}} \Lhier_\alpha.}$
\end{defn}

\begin{remark}
It is worth mentioning that the hierarchy as defined here does not look like the ``\emph{standard}'' one obtained through the syntactic approach. Firstly, $\defclosure(b)$ does not close under the fundamental operations, which is an approach we have chosen to take due to a lack of Strong Infinity in the background universe. Secondly, the \Lhier hierarchy does not stratify nicely by rank because $\defclosure(b)$ may potentially add sets that are not subsets of $b$, for example $x \times y$. Moreover, the ordinals of $\Lhier_\alpha$ may not be $\alpha$, for example we will have that $n \in \Lhier_{2n+1}$. These seem to be useful deficiencies in our weak context because they simplify some of the proofs and this is the approach taken by Barwise. However, we will address alternatives to this approach in \Cref{thm: with inf vs without} and \Cref{thm: comparison of defs}.
\end{remark}

\noindent We start the analysis by noting some of the basic properties of $\Lhier_\alpha$, which are all easy to see. The final property follows from the fact that each of the fundamental operations is $\Sigma$-definable and therefore $\defclosure$ is also a $\Sigma$-definable operation. The required statement then follows by $\Sigma$-recursion.

\begin{lemma}\label{thm: properties of L} $(\tf{IKP}^{\minusinf})$ For all ordinals $\alpha, \beta$: \begin{enumerate} \setlength \itemsep{0pt}
\item \label{L: property 1} If $\beta \in \alpha$ then $\Lhier_\beta \subseteq \Lhier_\alpha$,
\item\label{L: property 2} If $\beta \subseteq \alpha$ then $\Lhier_\beta \subseteq \Lhier_\alpha$,
\item\label{L: property 3} $\Lhier_\alpha \in \Lhier_{\alpha + 1}$,
\item If $x, y \in \Lhier_\alpha$ then for any $i \in I$, $\mathcal{F}_i(x, y) \in \Lhier_{\alpha + 1}$,
\item If for all $\beta \in \alpha$, $\beta + 1 \in \alpha$ then $\Lhier_\alpha$ is transitive,
\item \Lhier is transitive,
\item The operation $\alpha \mapsto \Lhier_\alpha$ is $\Sigma$-definable.
\end{enumerate}
\end{lemma}

\begin{theorem} \label{thm: L in IKP}
For every axiom $\varphi$ of $\tf{IKP}^{\minusinf}$, $\tf{IKP}^{\minusinf} \vdash \varphi^\Lhier$. Moreover,
\begin{center}
$\tf{IKP}^{\minusinf} + ``\emph{Strong Infinity}" \vdash (\emph{Strong Infinity})^\Lhier.$
\end{center}
\end{theorem}

\begin{proof}
The axioms of Extensionality and $\in$-induction follow from the fact that \Lhier is a transitive class. Pairing follows from $\funpair$ and Unions from $\fununion$. The Axiom of Empty Set holds because 
\[
\emptyset = \{ y \divline y \neq y \} = \Lhier_0 \in \Lhier_1.
\]
Bounded Separation follows from \Cref{thm: bounded sep}. Therefore it remains to prove Bounded Collection.

\indent Suppose that $\varphi(x, y,z)$ is a $\Sigma_0$-formula and, working in $\tf{IKP}^{\minusinf}$, assume that $a, z \in \Lhier$ and \mbox{$\forall x \in a \: \exists y \in \Lhier \: \varphi(x, y, z).$} Since $\varphi$ is $\Sigma_0$, it is absolute between $\Lhier$ and the background universe. Now, by definition of $\Lhier$, we get that 
\[
\forall x \in a \: \exists \beta \: \exists y \in \Lhier_\beta \: \varphi(x, y, z).
\]
Then, using the $\Sigma$-Collection principle in the background universe and the fact that the $\Lhier_\alpha$ hierarchy is constructed by a $\Sigma$-definable operation, there exists $\alpha$ such that 
\[
\forall x \in a \: \exists \beta \in \alpha \: \exists y \in \Lhier_\beta \: \varphi(x, y, z)
\]
which, by property \ref{L: property 1} of \Cref{thm: properties of L}, yields that $\forall x \in a \: \exists y \in \Lhier_\alpha \: \varphi(x, y, z).$ So, $b = \Lhier_\alpha$ is our required witness for Bounded Collection. \\

\noindent Now we work in the theory $\tf{IKP}^{\minusinf} +$ ``\emph{Strong Infinity}''. We first prove that $\omega \subseteq \Lhier_\omega$. This is done by induction by showing that for all $n \in \omega$, $n \in \Lhier_{2n + 1}$. Note that this holds because if $n \in \Lhier_{2n + 1}$ then 
\[
n + 1 = n \cup \{n\} = \fununion \big( n, \funpair(n, n) \big) \in \Lhier_{2n + 1 + 2}.
\]
Then $\omega = \{ n \in \Lhier_\omega \divline n = \emptyset \mathor \exists m \in n \: (n = m \cup \{m\}) \}$ will be in \Lhier by Bounded Separation.
\end{proof}

\noindent It is worth noting here that we will frequently claim that sets of the form 
\[
\{ z \in \Lhier_\delta \divline \varphi(u, z) \}
\]
are in $\Lhier_{\delta + k}$ for some $k \in \omega$, where $u \in \Lhier_\delta$ and $\varphi$ is a $\Sigma_0$-formula, without computing the required $k$. This $k$ could be computed by breaking down how $\varphi$ was built up using the fundamental operations, however this is often an unnecessarily tedious computation. Also, we will also often use formulae of the form $\exists i \in  \mathcal{I} \: \varphi(i)$ and these should be seen as an abbreviation for $\bigvee_{i \in \mathcal{I}} \varphi(i).$

An important property of the constructible universe is the viability of the Axiom of Constructibility; the axiom asserting that $\tf{V} = \Lhier$. We shall next prove that this axiom does indeed hold in \Lhier, that is:

\begin{theorem} \label{thm: constructibility in IKP}
$\tf{IKP}^{\minusinf} \vdash (\tf{V} = \Lhier)^\Lhier$.
\end{theorem}

\noindent Our method of proving this will closely follow the corresponding proof over \tf{IZF} by Lubarsky in \cite{LubarskyIntuitionisticL}. In particular, to prove the theorem it suffices to prove the following lemma which mirrors the version in \cite{LubarskyIntuitionisticL} using the syntactic version of the definability operator.

\begin{lemma}
\label{thm: approx in L}
$(\tf{IKP}^{\minusinf})$ For every ordinal $\alpha$ in \tf{V} there is an ordinal $\alpha^*$ in $\Lhier$ such that $\Lhier_\alpha = \Lhier_{\alpha^*}$.
\end{lemma}

\noindent In order to do this we define the operation of \emph{hereditary addition} on ordinals. This is necessary because in general it will not be true that $\beta \in \alpha$ implies that $\beta + 1 \in \alpha + 1$:

\begin{defn}[Lubarsky] \label{defn: hered addition} \index{Hereditary Addition} For ordinals $\alpha$ and $\gamma$, hereditary addition is defined inductively on $\alpha$ as 
\[
\alpha +_H \gamma \coloneqq \bigg( \bigcup \{ \beta +_H \gamma \divline \beta \in \alpha \} \cup \{ \alpha \} \bigg ) + \gamma,
\]
where ``$+$'' is the usual ordinal addition. We will also use the notation
\[
(\alpha +_H \gamma)^- \coloneqq \bigg( \bigcup \{ \beta +_H \gamma \divline \beta \in \alpha \} \cup \{ \alpha \} \bigg ).
\]
\end{defn}

\begin{proof}
First note that, by using the fundamental operations and \Cref{thm: bounded sep}, there is a fixed natural number $k$ such that for any ordinal $\tau$ and set $x \in \Lhier_\tau$, 
\[
\{ \gamma \in \Lhier_\tau \divline \gamma \in \tf{Ord} \mathand \defclosure(\Lhier_\gamma) \subseteq x \} \in \Lhier_{\tau + k}.
\]
We prove the lemma by induction on $\alpha$. So, suppose that $\forall \beta \in \alpha \; \exists \beta^* \in \Lhier (\Lhier_\beta = \Lhier_{\beta^*})$. Since the $\Lhier_\alpha$ hierarchy is constructed by a $\Sigma$-definable operation, by $\Sigma$-Collection, we can fix some ordinal $\delta$ in \tf{V} such that $\forall \beta \in \alpha \exists \beta^* \in \Lhier_\delta (\Lhier_\beta = \Lhier_{\beta^*})$.  Then, we can define $\alpha^*$ to be the ordinal
\[
\alpha^* \coloneqq \bigcup \{ \gamma \in \Lhier_{((\alpha \cup \delta) +_H k)^-} \divline \gamma \in \tf{Ord} \mathand \defclosure(\Lhier_\gamma) \subseteq \Lhier_{\alpha} \} \in \defclosure(\Lhier_{(\alpha \cup \delta) +_H k}).
\]
Note that we needed to take the union of the above set to ensure that $\alpha^*$ is an ordinal because we are not assuming $\Lhier_\gamma$ is a transitive set for arbitrary ordinals $\gamma$. Now clearly, for any ordinal $\alpha \in \tf{V}$, $\alpha^* \in \Lhier$. Moreover, if $x \in \Lhier_{\alpha^*}$ then there is some ordinal $\gamma \in \Lhier_{((\alpha \cup \delta) +_H k)^-}$ for which $\defclosure(\Lhier_\gamma) \subseteq \Lhier_{\alpha}$ and $x \in \defclosure(\Lhier_\gamma)$. Thus, $\Lhier_{\alpha^*} \subseteq \Lhier_\alpha$. For the reverse implication, we observe that by the definition of $\delta$, for any $\beta \in \alpha$ there is some $\beta^* \in \Lhier_\delta$ such that $\Lhier_\beta = \Lhier_{\beta^*}$. Thus $\beta^* \in \Lhier_{((\alpha \cup \delta) +_H k)^-}$ and $\defclosure(\Lhier_{\beta^*}) \subseteq \Lhier_\alpha$, so $\beta^* \in \alpha^*$. Hence,
\[
\Lhier_\alpha = \bigcup_{\beta \in \alpha} \defclosure(\Lhier_\beta) = \bigcup_{\beta \in \alpha} \defclosure(\Lhier_{\beta^*}) \subseteq \bigcup_{\gamma \in \alpha^*} \defclosure(\Lhier_\gamma) = \Lhier_{\alpha^*}.
\]

\end{proof} 

\noindent As mentioned at the beginning of this section, in the presence of Strong Infinity we can define the constructible universe using a different definability operator, $\Defop$, where $\Defop(b) \coloneqq \bigcup_{n \in \omega} \defclosure^n(b)$. This gives us an alternative way to construct the constructible universe, which we now show is equivalent as long as $\omega$ exists.

\begin{defn} For $\alpha$ an ordinal, $\Lhieromega_\alpha \coloneqq { \displaystyle \bigcup}_{\beta \in \alpha} \Defop(\Lhieromega_\beta)$ and 
\[
\Lhieromega \coloneqq \bigcup_{\alpha \in \tf{Ord}} \Lhieromega_\alpha.
\]
\end{defn}

\noindent As before, we can easily observe a few basic properties of this hierarchy:

\begin{prop} $(\tf{IKP}^{\minusinf} + ``\emph{Strong Infinity}")$ For all ordinals $\alpha, \beta$:
\begin{enumerate}
\item If $\beta \in \alpha$ then $\Lhieromega_\beta \subseteq \Lhieromega_\alpha$,
\item $\Lhieromega_\alpha \in \Lhieromega_{\alpha + 1}$,
\item $\Lhieromega_\alpha$ is transitive,
\item $\Lhieromega_\alpha$ is a model of Bounded Separation.
\end{enumerate}
\end{prop} 

\noindent It is possible to compare the hierarchies $\Lhier_\alpha$ and $\Lhieromega_\alpha$ via the following correspondence:

\begin{lemma}\label{thm: with inf vs without}
$(\tf{IKP}^{\minusinf} + ``\emph{Strong Infinity}")$ For any ordinal $\alpha$, $\Lhieromega_\alpha = \Lhier_{\omega \cdot \alpha}.$
\end{lemma}

\begin{proof}
We proceed by induction on $\alpha$. So assume that our claim holds for all $\beta \in \alpha$. First note that for any ordinal $\alpha$, $\omega \cdot \alpha = \{\omega \cdot \gamma + n \divline \gamma \in \alpha \mathand n \in \omega \}$ is an ordinal which is closed under successors. Then 

{
\setstretch{2} \begin{center} $ \displaystyle \begin{aligned} 
\Lhieromega_\alpha & = \bigcup_{\beta \in \alpha} \Defop(\Lhieromega_\beta) = \bigcup_{\beta \in \alpha} \bigcup_{n \in \omega} \defclosure^n(\Lhieromega_\beta) \\
& = \bigcup_{\beta \in \alpha} \bigcup_{n \in \omega} \defclosure^n(\Lhier_{\omega \cdot \beta}) = \bigcup_{\beta \in \alpha} \bigcup_{n \in \omega} \Lhier_{\omega \cdot \beta + n} \\
& = \Lhier_{\omega \cdot \alpha}. 
\end{aligned} $ \end{center} 
}

\end{proof}

\begin{cor} $(\tf{IKP}^{\minusinf} + ``\emph{Strong Infinity}")$ $\Lhieromega = \Lhier$.
\end{cor}

\noindent For completeness we briefly discuss how this method relates to the first, syntactic, approach we mentioned at the beginning of the chapter. We shall be sloppy in our presentation of the syntactic definability operator by using the ``$\models$'' symbol in our definition of definability instead of the more formal way this is presented in the previously mentioned references. We then defer to \cite{CrosillaPhD} for the formal way to do this in $\tf{IKP}^{\minusinf} +\emph{``Strong Infinity''}$. We remark here that the syntactic operator is the standard operator we shall use when taking the collection of definable subsets of a given set.

\begin{defn}
Say that a set $x$ is definable over a model $\langle \tf{M}, \in \rangle$ if there exists a formula $\varphi$ and $a_1, \dots, a_n \in \tf{M}$ such that 
\[
x = \{ y \in \tf{M} \divline \langle \tf{M}, \in \rangle \models \varphi[y, a_1, \dots, a_n] \}.
\]
We can then define the collection of definable subsets of \tf{M} as 
\[
\defop(\tf{M}) \coloneqq \{ x \subseteq \tf{M} \divline \textit{x is definable over} \: \langle \tf{M}, \in \rangle \}.
\]
\end{defn}

\noindent The constructible hierarchy can then be defined iteratively as 
\[
\mathbf{L}_\alpha \coloneqq { \displaystyle \bigcup}_{\beta \in \alpha} \defop(\mathbf{L}_\beta).
\]
Clearly, given $x, y \in \tf{M}$ and $i \in \mathcal{I}$, $\mathcal{F}_i(x, y)$ is definable from \tf{M} (even if it may take a few steps when $\mathcal{F}_i(x, y)$ is not a subset of \tf{M}, for example $x \times y \in \defop^3(\tf{M})$). Moreover, one can define the notion of being definable over $\langle \tf{M}, \in \rangle$ using only the fundamental operations, so the two universes they produce will be the same. To see the relationship between the two hierarchies, one can perform a careful analysis of the standard proof, for example Lemma VI$.1.17$ of \cite{DevlinBook}, which yields;

\begin{theorem}\label{thm: comparison of defs}
For every transitive set \tf{M}: 
\begin{eqnarray*} 
\textnormal{def}(\tf{M}) & = & \textnormal{Def}(\tf{M}) \cap \mathcal{P}(\tf{M}) \\
& = & \bigcup_{n \in \omega} \defclosure^n(\tf{M}) \cap \mathcal{P}(\tf{M}).
\end{eqnarray*}
\end{theorem}

\noindent Therefore, in the theory $\tf{IKP}^{\minusinf} + ``\emph{Strong Infinity}"$ the two standard formulations of the constructible hierarchy are equivalent. One of the main benefits for using the formulation in terms of the fundamental operations is to avoid the use of Strong Infinity in the construction. The second benefit for using our formulation is because of the versatility of these operations over \tf{ZF}. A notable example is that it allows us to define when an inner model satisfies \tf{ZF}. This occurs when the inner model is closed under fundamental operations and it satisfies a property known as \emph{almost universality}. We shall see in the next section that an analogous result holds in \tf{IZF}.

\pagebreak
\section{External Cumulative Hierarchies}

In this section we shall show that if  \tf{V} satisfies \tf{IZF} then so does \Lhieromega. This could be done by a very similar repetition of the analysis in the \tf{IKP} case however we will take a different approach here in order to derive further axiomatic properties under \tf{IZF}.

The main theorem of this section is adapted from Theorem 13.9 of \cite{JechSetTheory}. The essence of that theorem is that if \tf{V} is a model of \tf{ZF} and \tf{M} contains all of the ordinals then \tf{M} being a model of \tf{ZF} can be expressed by a single first-order sentence. On the face of it, the theorem we present here will be slightly weaker than this because it will require the additional assumption that \tf{M} has an \emph{external cumulative hierarchy}. 

\begin{defn} \index{External Cumulative Hierarchy} Let $\tf{M} \subseteq \tf{N}$. We say that \tf{M} has an \emph{external cumulative hierarchy} (e.c.h.) in \tf{N} if there exists a sequence of sets $\langle \tf{M}_\alpha \divline \alpha \in \tf{Ord} \cap \tf{N} \rangle$, definable in \tf{N}, such that: \begin{itemize}
\item For every $\alpha \in \tf{Ord} \cap \tf{N}$, $\tf{M}_\alpha \in \tf{M}$,
\item ${\displaystyle \tf{M} = \bigcup  \{ \tf{M}_\alpha \divline \alpha \in \tf{Ord} \: \cap \: \tf{N} \} }$,
\item If $\beta \in \alpha$ then $\tf{M}_\beta \subseteq \tf{M}_\alpha$.
\end{itemize}
\end{defn}

\noindent We say that \tf{M} has an e.c.h. when $\tf{N} = \tf{V}$. It is worth remarking that if \tf{M} is an inner model of \tf{IZF}, that is a model of \tf{IZF} containing all of the ordinals, then \tf{M} will have an external cumulative hierarchy given by the standard rank hierarchy which can be defined as follows:

\vbox{
\begin{defn} Define the rank of $a$, $\ff{rank}(a)$, recursively as follows: 
\[
\ff{rank}(a) \coloneqq \bigcup \{ \ff{rank}(x) + 1 \divline x \in a \},
\]
where $z + 1 \coloneqq z \cup \{z\}$.
\end{defn}
} 

\noindent Note that one can easily prove that for any set $a$, $\ff{rank}(a)$ is an ordinal and for any ordinal $\alpha$, $\ff{rank}(\alpha) = \alpha.$ 

\begin{defn} 
For $\alpha$ an ordinal, $ \displaystyle \tf{V}_\alpha \coloneqq \bigcup_{\beta \in \alpha} \mathcal{P}(\tf{V}_\beta).$ 
\end{defn}

\begin{prop}\label{ranks} 
For any set $a$, $a \subseteq \tf{V}_{\ff{rank}(a)}.$ 
\end{prop} 

\begin{proof} This is formally proved by induction on rank noting that for any $x \in a$, if $x \subseteq \tf{V}_{\ff{rank}(x)}$ then $x \in \mathcal{P}(\tf{V}_{\ff{rank}(x)}) \subseteq \tf{V}_{\ff{rank}(a)}$. 
\end{proof} 

\noindent Therefore, if \tf{M} is a definable class which is an inner model of \tf{IZF}, then $\langle \tf{V}_\alpha ^{\tf{M}} \divline \alpha \in \tf{Ord} \rangle$ defines an e.c.h. which is moreover uniformly definable. Also, by construction, we have that $\langle \Lhieromega_\alpha \divline \alpha \in \tf{Ord} \rangle$ is an e.c.h. for \Lhieromega even though it is unlikely to be the case that \Lhieromega contains all of the ordinals of \tf{V}. 

\begin{defn} \index{Almost Universal} Let $\tf{M} \subseteq \tf{N}$. We say that \tf{M} is \emph{almost universal} in \tf{N} if for any $x \in \tf{N}$, if $x \subseteq \tf{M}$ then there exists some $y \in \tf{M}$ such that $x \subseteq y$.
\end{defn}

\noindent Furthermore, it is worth pointing out that classically, for almost universal models, having an e.c.h. implies that the two models have the same ordinals. This will be stated in \tf{ZF} for simplicity, but in reality only requires very basic set theory and some amount of separation and replacement with regards to the hierarchy.

\begin{prop} Suppose that $\tf{M} \subseteq \tf{N}$ are transitive models of \tf{ZF} and \tf{M} is almost universal in \tf{N}. If \tf{M} has an e.c.h. in \tf{N} then $\tf{Ord} \cap \tf{M} = \tf{Ord} \cap \tf{N}$.
\end{prop}

\begin{proof}
Let $\langle \tf{M}_\alpha \divline \alpha \in \tf{Ord} \cap \tf{N} \rangle$ be an external cumulative hierarchy. We shall prove inductively that for any ordinal $\gamma \in \tf{N}$ there is an ordinal $\beta \in \tf{N}$ such that $\gamma \subseteq \tf{M}_\beta$. Then, since $\tf{M}_\beta \in \tf{M}$, either $\gamma = \tf{M}_\beta \cap \tf{Ord}$ or $\gamma \in \tf{M}_\beta \cap \tf{Ord}$. 
Almost universality will then allow us to take some transitive set $y \in \tf{M}$ covering the set $\tf{M}_\beta \cap \tf{Ord}$ and then Bounded Separation in \tf{M} yields that $\gamma \in \tf{Ord} \cap \tf{M}$.

Working in \tf{N}, to prove the claim first note that, by induction,
\[
\forall \alpha \in \gamma \: \exists \tau_\alpha \in \tf{N} \: \alpha \in \tf{M}_{\tau_\alpha}.
\]
So, by Collection and the assumption that the hierarchy is cumulative, there is some ordinal $\beta$ such that for all $\alpha \in \gamma$, $\alpha \in \tf{M}_\beta$ and therefore $\gamma \subseteq \tf{M}_\beta$. 
\end{proof}

\begin{remark}
For clarity, the point where we used excluded middle in the above proof was when we asserted that $\gamma \subseteq \tf{M}_\beta$ implies that $\gamma \in \tf{M}$. If $\gamma = \tf{M}_\beta \cap \tf{Ord}$ then $\gamma$ will be in $\tf{Ord} \cap \tf{M}$ as 
\[
\gamma = \{ \alpha \in \tf{M}_\beta \divline \alpha \textit{ is a transitive set of transitive sets} \}.
\]
However, if $\gamma$ is a proper subset of $\tf{M}_\beta \cap \tf{Ord}$ then we require linearity of the ordinals to conclude that $\gamma$ is an element of this set. To see an example of this, consider an arbitrary truth value $x \subseteq 1$. It is quite plausible for such an ordinal to be a proper subset of some $\tf{M}_\beta \cap \tf{Ord}$ since this set could contain $1$ but there is no reason to assume that $x$ itself is a member of this set.  
\end{remark}

\noindent We now present our adaptation of the theorem from Jech. The assumption of an e.c.h. is not necessary for the right-to-left implication of this theorem but it is needed in order to prove that our model \tf{M} is almost universal in \tf{V}. This is because, while we can use the rank hierarchy of \tf{M} to show that if a set $a \in \tf{V}$ is a subset of \tf{M} there is some $\beta \in \tf{Ord}$ such that $a \subseteq {\displaystyle \bigcup}_{\alpha \in \beta} \tf{V}_\alpha ^{\tf{M}},$ it does not seem possible to show that this union is in fact a set in \tf{M} because there is no reason why $\beta$ should be in \tf{M}. 

However, having an e.c.h. seems to be a reasonable additional assumption since in most cases our model \tf{M} will be built up iteratively over the ordinals in \tf{V}, which gives a very natural hierarchy. Notably, we have that $\Lhieromega \coloneqq {\displaystyle \bigcup}_{\alpha \in \tf{Ord} \cap \tf{V}} \Lhieromega_\alpha.$

\ConditionforMtoModelIZF

\begin{proof}
For the left-to-right implication, first see that, if \tf{M} is a model of \tf{IZF}, then \tf{M} is certainly closed under the fundamental operations because they are $\Sigma$-definable. For almost universality, let $\langle \tf{M}_\alpha \divline \alpha \in \tf{Ord} \rangle$ be an e.c.h. and suppose that $a \in \tf{V}$ with $a \subseteq \tf{M}$. Then we have that 
\[
\forall x \in a \: \exists \alpha \in \tf{Ord} \: (x \in \tf{M}_\alpha).
\]
So, by Collection in \tf{V}, there is some set $b \in \tf{V}$ such that 
\[
\forall x \in a \: \exists \alpha \in b \: (x \in \tf{M}_\alpha).
\]
Taking $\beta = \ff{trcl}(b) \cap \tf{Ord}$ we can then conclude that $a \subseteq \tf{M}_\beta$, and $\tf{M}_\beta$ is a set in \tf{M} by the assumption that the $\tf{M}_\alpha$'s form an external cumulative hierarchy. 

For the reverse implication, since \tf{M} is transitive it is a model of Extensionality and $\in$-induction. Also, Pairing and Unions follow from \tf{M} being closed under $\funpair$ and $\fununion$ while the Axiom of Infinity will follow from the proof of \Cref{thm: L in IKP} along with an instance of Bounded Separation. We now proceed to prove the other axioms. 

\smallskip

\noindent \textbf{Bounded Separation:} This follows from the same argument as given in \Cref{thm: bounded sep} because any $\Sigma_0$-formula can be expressed using the fundamental operations. 

\smallskip

\noindent \textbf{Power Set:} Let $a$ be a set and note that $\mathcal{P}(a) \cap \tf{M}$ is a set in \tf{V}, and a subset of \tf{M}. So, by almost universality, we can fix some $b \in \tf{M}$ such that $\mathcal{P}(a) \cap \tf{M} \subseteq b$. But then 
\[
\mathcal{P}^{\tf{M}}(a) = \{ x \in b \divline x \subseteq a \},
\]
which is a set in \tf{M} by Bounded Separation in \tf{M}.  \\

\noindent \textbf{Collection:} Let $a$ be in \tf{M} and suppose that, in \tf{M}, $\forall x \in a \: \exists y \in \tf{M} \: \varphi(x, y, u).$ By Collection in \tf{V} we can fix some set $b'$ such that 
\[
\tf{V} \models \forall x \in a \: \exists y \in b' \: \varphi^{\tf{M}}(x, y, u)
\]
and, by almost universality, we can fix $b \in \tf{M}$ such that $b \supseteq b'$ yielding, in \tf{M}, 
\[
\tf{M} \models \forall x \in a \: \exists y \in b \: \varphi(x, y, u).
\]

\noindent \textbf{Separation:} This is shown by induction on the complexity of the formula $\varphi$. $\Sigma_0$-formulae and all cases except for those involving quantifiers follow immediately from the consequences of the fundamental operations. 

So suppose that $\varphi(x, u) \equiv \exists v \: \psi(x, u, v)$. Using Separation in \tf{V}, define $a'$ to be 
\[
a' \coloneqq \{ x \in a \divline \varphi^{\tf{M}} (x, u) \}.
\]

\noindent Then, 
\[
\tf{V} \models \forall x \in a' \: \exists v \in \tf{M} \: \psi^{\tf{M}}(x, u, v).
\]
So, by Collection in \tf{V}, there exists some set $b' \subseteq \tf{M}$ such that
\[
\tf{V} \models \forall x \in a' \: \exists v \in b' \: \psi^{\tf{M}}(x, u, v).
\]
By almost universality, take $b \in \tf{M}$ such that $b' \subseteq b$. Then 
\[
\tf{V} \models \forall x \in a' \: \exists v \in b \: \psi^{\tf{M}}(x, u, v).
\]
Now, by the inductive hypothesis, we have that 
\[
y \coloneqq \{ \langle x, v \rangle \in a \times b \divline \psi^{\tf{M}} (x, u, v) \} \in \tf{M},
\]
and thus, using $\fundom$, 
\[
z = \ff{dom}(y) = \{ x \in a \divline \exists v \in b \: \psi^{\tf{M}}(x, u, v) \} \in \tf{M}.
\]

\bigskip
\bigskip

\noindent For the final case, suppose that $\varphi(x, u) \equiv \forall v \: \psi(x, u, v)$. For $r \in \tf{M}$ let 
\[
y_r \coloneqq \{ x \in a \divline \forall v \in r \: \psi^{\tf{M}} (x, u, v) \} \in \tf{V}.
\]
Using the inductive hypothesis and the proof of case \ref{thm: for all x in b is t-form} of \Cref{lemma:BarwisetFormula} we shall show that $y_r \in \tf{M}$ for every $r \in \tf{M}$. Firstly, by the inductive hypothesis, $c \coloneqq \{ \langle z, y, x \rangle \in r \times \{u\} \times a \divline \psi^{\tf{M}}(x, y, z) \}$ is in $\tf{M}$. Next, by definition, $\funforall(c, r) = \{ \{ \langle u, x \rangle \divline \langle z, u, x \rangle \in c \} \divline z \in r \}$. Thus, $\funint( \{u\} \times a, \funforall(c, r) ) =$ $\{ \langle u, x \rangle \divline x \in a \wedge \forall z \in r \; \psi^{\tf{M}}(x, u, z)\}$ from which it follows that $y_r = \funran(\funint(\{u\} \times a, \funforall(c, r)))$ is in \tf{M} since \tf{M} is closed under the fundamental operations.

Furthermore, it is obvious that if $s \subseteq r$ then $y_r \subseteq y_s$. We aim to show that $y_\tf{M}$, which is defined in the same way by taking a universal quantifier over \tf{M}, is in \tf{M} by showing that it is equal to $y_r$ for some $r \in \tf{M}$. To do this, we begin by defining 
\[
Y \coloneqq \{ z \in \mathcal{P}(a) \divline \exists r \in \tf{M} \: (z = y_r) \}.
\]
So 
\[
\tf{V} \models \forall z \in Y \: \exists r \in \tf{M} \: (z = y_r).
\]
Therefore, by Collection in \tf{V}, there is some set $d' \subseteq \tf{M}$ such that 
\[
\tf{V} \models \forall z \in Y \: \exists r \in d' \: (z = y_r).
\]
Taking the transitive closure if necessary, by almost universality fix a transitive set $d \in \tf{M}$ such that $d' \subseteq d$. We claim that $y_d = y_\tf{M}$. Firstly, since $d \subseteq \tf{M}$, $y_\tf{M} \subseteq y_d$. For the reverse direction, let $x \in y_d$ and let $v \in \tf{M}$. Then $y_{\{v\}} \in Y$ so we can fix $r \in d$ such that $y_{\{v\}} = y_r$. Therefore, since $r \subseteq d$, $y_d \subseteq y_r = y_{\{v\}}$ so $x \in y_{\{v\}}$ and, by construction, $\psi^{\tf{M}}(x, v, u)$. Finally, since $v$ was arbitrary, $x \in y_\tf{M}$ as required.
\end{proof}

\begin{cor}
For every axiom $\varphi$ of $\tf{IZF}$, $\tf{IZF} \vdash \varphi^{\Lhieromega}$.
\end{cor}

\section{Exponentiation in the Constructible Universe} \label{section:ExpinL}
In order to show that $\tf{CZF}_\mathcal{P}$ cannot prove that the Axiom of Exponentiation holds in \tf{L} we shall make use of a specific realizability model. This model is based upon general set recursive functions where a realizer for an existential statement will provide a set of witnesses for the existential quantifier, rather than just a single witness. When combined with truth (for implications) this will then produce an effective procedure for constructing what are known as $E_\wp$-recursive functions which witness sentences provable in $\tf{CZF}_\mathcal{P}$. This construction was originally derived in \cite{RathjenExistenceProperty} and we refer to that paper for all necessary details. 

In order to give the proof it is necessary to first briefly outline the notions of $E$-recursive and $E_\wp$-recursive functions and the realizability construction that we shall use. This shall be done over the next two subsections.

\subsection{\texorpdfstring{$E$}{E}-recursive functions}

$E$-recursion is a way to extend Kleene's notion of computability, $\{e\}(x)$, to apply to arbitrary sets. This generalised notion of recursion, also known as set recursion, is also discussed in detail in sources such as \cite{NormannSetRecursion} or chapter 10 of \cite{SacksRecursionTheory}. The idea is that, for arbitrary sets $a$ and $x$, $a$ should be thought of as some form of Turing Machine which takes as input $x$. In order to avoid any possible confusion with the singleton set, we shall instead use $[a](c)$ for the result of this computation. $E_\wp$-recursion will be an extension of this generalised notion to also treat power sets as an initial function. This latter notion has also been previously studied by authors such as Moschovakis \cite{MoschovakisRecursion} and Moss \cite{MossPowerSetRecursion}.

Here we shall use Kleene's notion of \emph{strong equality}, $\simeq$, which is that for two partial functions $f$ and $g$, $f(x) \simeq g(x)$ if and only if neither $f(x)$ nor $g(x)$ are defined, or they are both defined and equal.

\begin{defn}[\cite{RathjenExistenceProperty} 2.9] (\tf{IKP}) \label{definition:SpecialIndices}
First, we select distinct non-zero natural numbers \textbf{k}, \textbf{s}, \textbf{p}, \textbf{p}\textsubscript{\textbf{0}}, \textbf{p}\textsubscript{\textbf{1}}, \textbf{s}\textsubscript{\textbf{N}}, \textbf{p}\textsubscript{\textbf{N}}, \textbf{d}\textsubscript{\textbf{N}}, $\boldsymbol{\bar{0}}$, $\boldsymbol{\bar{\omega}}$, $\boldsymbol{\pi}$, $\boldsymbol{\nu}$, $\boldsymbol{\gamma}$, $\boldsymbol{\rho}$, \textbf{i}\textsubscript{\textbf{1}}, \textbf{i}\textsubscript{\textbf{2}} and \textbf{i}\textsubscript{\textbf{3}} which will provide indices for special $E$-recursive class functions. Let $\mathbb{N} \coloneqq \omega$.

We shall inductively define a class $\mathbb{E}$ of triples $\langle e, x, y \rangle$. Instead of saying ``$\langle e, x, y \rangle \in \mathbb{E}$'' we shall write ``$[e](x) \simeq y$''. Next, we shall use $[e](x_1, \dots x_n) \simeq y$ to indicate that 
\[
[e](x_1) \simeq \langle e, x_1 \rangle \: \mathand \: [\langle e, x_1 \rangle](x_2) \simeq \langle e, x_1, x_2 \: \rangle \: \mathand \dots \: \mathand \: [\langle e, x_1, \dots, x_{n-1} \rangle](x_n) \simeq y.
\]
Finally, we shall say that $[e](x)$ is defined, written $[e](x)\downwards$, iff $[e](x) \simeq y$ for some $y$. The relation $[e](x) \simeq y$ is then specified by the following clauses: 
\begin{align*}
\intertext{Indices for applicative structures (APP):}
\begin{split}
[\textbf{k}](x,y) & \; \simeq \; x \\
[\textbf{p}](x,y) & \; \simeq \; \langle x, y \rangle \\ 
[\textbf{p}_{\textbf{N}}](0) & \; \simeq \; 0 \\
[\textbf{p}_{\textbf{0}}](x) & \; \simeq \; 1^{st}(x) \\
[\textbf{d}_{\textbf{N}}](n, m, x, y) & \; \simeq \; x \textnormal{ if } n, m \in \mathbb{N} \textnormal{ and } n = m \\
[\boldsymbol{\bar{0}}](x) & \; \simeq \; 0 \\
\end{split}
\begin{split}
[\textbf{s}](x,y,z) & \; \simeq \; [[x](z)]([y](z)) \\
[\textbf{s}_{\textbf{N}}](n) & \; \simeq \; n + 1 \textnormal{ if } n \in \mathbb{N} \\
[\textbf{p}_{\textbf{N}}](n+1) & \; \simeq \; n \textnormal{ if } n \in \mathbb{N} \\
[\textbf{p}_{\textbf{1}}](x) & \; \simeq \; 2^{nd}(x) \\
[\textbf{d}_{\textbf{N}}](n, m, x, y) & \; \simeq \; y \textnormal{ if } n, m \in \mathbb{N} \textnormal{ and } n \neq m \\
[\boldsymbol{\bar{\omega}}](x) & \; \simeq \; \omega \\
\end{split}
\intertext{Indices for set-theoretic axioms:}
\begin{split}
[\boldsymbol{\pi}](x, y) & \; \simeq \; \{x, y\} \\
[\boldsymbol{\gamma}](x, y) & \; \simeq \; x \cap \bigcap y \\
\end{split}
\begin{split}
[\boldsymbol{\nu}](x) & \; \simeq \; \bigcup x \\
[\boldsymbol{\rho}](x, y) & \; \simeq \; \{ [x](u) \divline u \in y \} \\
& \phantom{\; \simeq \;} \textnormal{ if } [x](u) \textnormal{ is defined for all } u \in y \\
\end{split}
\intertext{Indices for equality axioms:}
\begin{split}
[\textbf{i}_{\textbf{1}}](x, y, z) & \; \simeq \; \{ u \in x \divline y \in z \} \\
[\textbf{i}_{\textbf{3}}](x, y, z) & \; \simeq \; \{ u \in x \divline u \in y \rightarrow z \in u \} \\
\end{split}
\begin{split}
[\textbf{i}_{\textbf{2}}](x, y, z) & \; \simeq \; \{ u \in x \divline u \in y \rightarrow u \in z \} \\
{}&{} \\
\end{split}
\end{align*}

\end{defn}

\begin{defn}[\cite{RathjenExistenceProperty} 2.11]
\emph{Application terms} are defined inductively as:
\begin{enumerate}
\item The distinguished indices indicated in \Cref{definition:SpecialIndices} are all application terms,
\item variables are application terms,
\item if $s$ and $t$ are application terms then so is $[s](t)$.
\end{enumerate}
An application term is then called \emph{closed} if it does not contain any variables.
\end{defn}

\begin{defn}[\cite{RathjenExistenceProperty} 2.13]
A partial $n$-place class function $\Upsilon$ is called an $E$\emph{-recursive partial function} if there exists a closed application term $t_{\Upsilon}$ such that
\[
\ff{dom}(\Upsilon) = \{ (a_1, \dots, a_n) \divline [t_{\Upsilon}](a_1, \dots, a_n)\downarrow \}
\]
and for all sets $(a_1, \dots, a_n) \in \ff{dom}(\Upsilon)$,
\[
[t_{\Upsilon}](a_1, \dots, a_n) \; \simeq \; \Upsilon(a_1, \dots, a_n).
\]
When this holds, $t_\Upsilon$ is said to be an \emph{index} for $\Upsilon$ and we will freely confuse the two.
\end{defn}

\noindent When working with $\tf{IKP}(\mathcal{P})$ we can extend the above notion of recursive functions to $E_{\wp}$\emph{-recursive functions}, which will make the function $\mathcal{P}(x) = \{u \divline u \subseteq x \}$ ``\emph{computable}''. This is done by specifying an additional non-zero natural number $\bar{\boldsymbol\wp}$ and adding the clause
\[
[\bar{\boldsymbol\wp}](x) \; \simeq \; \mathcal{P}(x)
\]
to \Cref{definition:SpecialIndices}, giving rise to the class $\mathbb{E}_{\wp}$. The notions of $\wp$-application terms and closed $\wp$-application terms are then defined analogously.

Finally, we can state the main technical result for these functions, which is Lemmas 2.20 and 2.21 of \cite{RathjenExistenceProperty}.

\begin{lemma}[\cite{RathjenExistenceProperty} 2.20, 2.21] \,
\begin{enumerate}
\item For each $\Sigma_0$-formula $\varphi(a, y)$ there is a closed application term $t_\varphi$ such that
\[
\tf{IKP} \vdash [t_\varphi](a, b)\downwards \mathand \forall u \big( u \in [t_\varphi](a, b) \longleftrightarrow (u \in b \mathand \varphi(a, u)) \big).
\]
\item For each $\Sigma^{\mathcal{P}}_0$-formula $\varphi(a, y)$ there is a closed $\wp$-application term $t_\varphi$ such that
\[
\tf{IKP}(\mathcal{P}) \vdash [t_\varphi](a, b)\downwards \mathand \forall u \big( u \in [t_\varphi](a, b) \longleftrightarrow (u \in b \mathand \varphi(a, u)) \big).
\]\end{enumerate}
\end{lemma}

\subsection{Realizability Preliminaries}

We can now give the definition of the notion of realizability that we shall use, $\rlzt$, which will be a relation between sets and formulae. Following the notation of the second author in \cite{RathjenExistenceProperty} we shall use the subscript $\mathfrak{w}$ to indicate \emph{weakening} the realizer for an existential statement to only be a set of witnesses for the existential quantifier and the subscript $\mathfrak{t}$ to denote the addition of \emph{truth} to the realizability interpretation. In this relation, we will treat bounded quantifiers as quantifiers in their own right. In the definition we shall use the following notation:

\begin{itemize}
\item $\pairl e$ and $\pairr e$ for $1^{st}(e)$ and $2^{nd}(e)$ respectively (instead of $[\pairl](e)$ and $[\pairr](e)$),
\item $[a](b) \rlzt \varphi$ for $\exists x \big( [a](b) \simeq x \mathand x \rlzt \varphi \big)$, where $[a](b) \simeq x$ denotes $\langle a, b, x \rangle \in \mathbb{E}$.
\end{itemize}  

\begin{defn}[\cite{RathjenExistenceProperty} 3.1] \label{definition:realizabilitywt} \,
{\begin{align*}
a & \rlzt \varphi && \textnormal{iff \quad} \varphi \textnormal{ is true, whenever } \varphi \textnormal{ is an atomic formula},\\
a & \rlzt \varphi \mathand \psi && \textnormal{iff \quad} \pairl a \rlzt \varphi \mathand \pairr a \rlzt \psi, \\
a & \rlzt \varphi \mathor \psi && \textnormal{iff \quad} \exists u ( u \in a ) \mathand \forall d \in a \: \big( ( \pairl d = 0 \mathand \pairr d \rlzt \varphi ) \mathor (\pairl d = 1 \mathand \pairr d \rlzt \psi) \big),  \\
a & \rlzt \varphi \rightarrow \psi && \textnormal{iff \quad} (\varphi \rightarrow \psi) \mathand \forall c ( c \rlzt \varphi \rightarrow [a](c) \rlzt \psi), \\
a & \rlzt \forall x \in b \: \varphi(x) && \textnormal{iff \quad}  \forall c \in b \: ([a](c) \rlzt \varphi[x/c]), \\
a & \rlzt \exists x \in b \: \varphi(x) && \textnormal{iff \quad} \exists u ( u \in a ) \mathand \forall d \in a \: ( \pairl d \in b \mathand \pairr d \rlzt \varphi[x/\pairl d] ), \\
a & \rlzt \forall x \: \varphi(x) && \textnormal{iff \quad}  \forall c \: [a](c) \rlzt \varphi[x/c], \\
a & \rlzt \exists x \: \varphi(x) && \textnormal{iff \quad} \exists u ( u \in a ) \mathand \forall d \in a (  \pairr d \rlzt \varphi[x/\pairl d] ), \\
& \rlzt \varphi && \textnormal{iff \quad} \exists a \: a \rlzt \varphi.
\end{align*}}
\end{defn}

\begin{terminology}
In \cite{RathjenExistenceProperty} the abbreviation ``$a \neq \emptyset$'' is used in Definition 3.1 to express the positive statement $\exists u (u \in a)$. We have chosen to write this statement explicitly in order to avoid any confusion.
\end{terminology}

\noindent Note that the above realizability notion was based on $E$-computability, and this will be used to provide a translation of $\tf{CZF}^{\uminus}$. In order to provide a translation of $\tf{CZF}_\mathcal{P}$ we will instead use the corresponding realizability notions based on $E_\wp$-computability.

\begin{defn} \label{definition:realizability_p}
$a \rlztp \varphi$ is defined recursively by the same recursion as in \Cref{definition:realizabilitywt} but using $E_\wp$-computability where $[a](b) \simeq c$ now stands for $\langle a, b, c \rangle \in \mathbb{E}_\wp$.
\end{defn}

\noindent By induction on the construction of formulae, one can show that this notion of realizability implies truth over both of our constructive theories.

\begin{prop}[\cite{RathjenExistenceProperty} 3.3] \label{theorem:RealizeImpliesTruth}
Let $\varphi(x_1, \dots, x_n)$ be a formula with free variables among $x_1, \dots, x_n$. Then
\begin{enumerate}
\item $\tf{CZF}^{\uminus} \vdash (\exists e \: e \rlzt \varphi(x_1, \dots, x_n)) \rightarrow \varphi(x_1, \dots, x_n).$
\item \label{theorem:RealizeImpliesTruth3} $\tf{CZF}_{\mathcal{P}} \vdash (\exists e \: e \rlztp \varphi(x_1, \dots, x_n)) \rightarrow \varphi(x_1, \dots, x_n).$
\end{enumerate}
\end{prop}

\noindent Next, by induction on the complexity of formulae, one can show that within our constructive theories our realizability structure realizes every statement provable from the given theory. This is done by first ensuring that the interpretation preserves logical inferences and then that we are able to realize every axiom of $\tf{CZF}^{\uminus}$ ($\tf{CZF}_{\mathcal{P}}$) via an $E$-recursive ($E_\wp$-recursive) function. 

\begin{theorem}[\cite{RathjenExistenceProperty} 3.7, 3.9] \label{theorem:EffectivelyConstructIndex}
Let $\varphi(x_1, \dots, x_n)$ be a formula with free variables among $x_1, \dots, x_n$.
\begin{enumerate}
\item If $\tf{CZF}^{\uminus} \vdash \varphi(x_1, \dots, x_n)$ then one can effectively construct an index of an $E$-recursive function $f$ such that $\tf{CZF}^{\uminus} \vdash \forall a_1 \dots \forall a_n \: f(a_1, \dots a_n) \rlzt \varphi(a_1, \dots, a_n)$.
\item \label{theorem:EffectivelyConstructIndex3} If $\tf{CZF}_{\mathcal{P}} \vdash \varphi(x_1, \dots, x_n)$ then one can effectively construct an index of an $E_\wp$-recursive function $f$ such that $\tf{CZF}_{\mathcal{P}} \vdash \forall a_1 \dots \forall a_n \: f(a_1, \dots a_n) \rlztp \varphi(a_1, \dots, a_n)$.
\end{enumerate}  
\end{theorem}

\noindent Moreover, if one removes the ``\textit{truth}'' requirement, that is by defining the interpretation $\rlz$ using the same conditions as in \Cref{definition:realizabilitywt} except that the statement for implication is replaced by
\[
a \rlz \varphi \rightarrow \psi \textnormal{\quad iff \quad} \forall c ( c \rlz \varphi \rightarrow [a](c) \rlz \psi),
\]

\noindent then one can perform a similar analysis to give a realizability interpretation of $\tf{CZF}^{\uminus}$ and $\tf{CZF}_\mathcal{P}$ over the background theories of \tf{IKP} and $\tf{IKP}(\mathcal{P})$. Moreover, this interpretation will imply truth for bounded formulae. These translations can then be used to show that $\tf{CZF}^{\uminus}$ proves the same $\Pi_2$ assertions as $\tf{IKP}$. Note that a formula is said to be $\Pi_2$ if it is of the form $\forall x \: \exists y \: \varphi(x, y, z)$ where $\varphi(x, y, z)$ is $\Sigma_0$. The class $\Pi^{\mathcal{P}}_2$ is defined analogously by allowing  $\varphi(x, y, z)$ to be of $\Sigma_0^{\mathcal{P}}$ form. 

\begin{theorem}[\cite{RathjenExistenceProperty} 4.8] \label{theorem:CZFPisPi2Conservative} \,
\begin{enumerate}
\item $\tf{CZF}^{\uminus}$ is conservative over $\tf{IKP}$ for $\Pi_2$-sentences.
\item $\tf{CZF}_\mathcal{P}$ is conservative over $\tf{IKP}(\mathcal{P})$ for $\Pi^{\mathcal{P}}_2$-sentences.
\end{enumerate}
\end{theorem}

\noindent Via an ordinal analysis using infinitary cut-free derivations, a strong form of reflection for $\Sigma^\mathcal{P}$-formulae is obtained in Theorem 8.2 of \cite{RathjenRelativizedOrdinalAnalysis}. Namely, from any proof of a $\Sigma^\mathcal{P}$-sentence $\varphi$ one can effectively extract an ordinal (representation) $\sigma$ below the Bachmann-Howard ordinal, $\psi_{\Omega}(\varepsilon_{\Omega+1})$, such that $\tf{KP}(\mathcal{P})$ proves that the von Neumann hierarchy exists up to level $\sigma$ and $\varphi$ holds in $\tf{V}_\sigma$. It was further shown in \cite{Cook-Rathjen} that this technique can be transfered to the intuitionistic counterpart of $\tf{KP}(\mathcal{P})$, $\tf{IKP}(\mathcal{P})$.

\begin{theorem}[\cite{Cook-Rathjen} 3.27] \label{theorem:SigmaReflectstoBH}
Let $\varphi$ be a $\Sigma^\mathcal{P}$-sentence. If $\tf{IKP}(\mathcal{P}) \vdash \varphi$ then there exists some ordinal term $\sigma < \psi_{\Omega}(\varepsilon_{\Omega+1})$, which can be computed from the derivation, such that $\tf{V}_\sigma \models \varphi$.
\end{theorem}

\subsection{Exponentiation in \texorpdfstring{$\tf{CZF}(\mathcal{P})$}{CZF(P)}}

Using this machinery we are now able to show that $\tf{CZF}_\mathcal{P}$ is unable to prove that the Axiom of Exponentiation holds in \tf{L}. In particular, this means that \tf{CZF} is unable to show that \tf{L} satisfies every axiom of \tf{CZF}. 

\ExpFailsinLofCZF

\begin{proof}
It will suffice to show that $\tf{CZF}_\mathcal{P}$ does not prove the following: \textit{The statement ``the collection of all functions from $\omega$ to $\omega$ forms a set'' relativizes to \tf{L}}. Aiming for a contradiction, suppose that $\tf{CZF}_\mathcal{P}$ was able to prove this statement.

Let $\varphi_{\ff{nat}}(u)$ be a $\Sigma_0$ formula which uniquely defines $\omega$, that is to say $\tf{IKP} \vdash \exists ! u \varphi_{\ff{nat}}(u)$. For example:
\[
\varphi_{\ff{nat}}(u) \equiv \exists x \in u \: \forall y \in x (y \neq y) \mathand \forall x \in u (x \cup \{ x\} \in u) \mathand \forall x \in u (x = \emptyset \mathor \exists y \in u (x = y \cup \{y\})).
\]
Now, by our assumption,
\[
\tf{CZF}_\mathcal{P} \vdash \exists u, z \in \tf{L} \Big( \varphi_{\ff{nat}}(u) \mathand \forall f \in \prescript{u}{}u \big( f \in \tf{L} \rightarrow f \in z \big) \Big).
\]
Therefore,
\[
\tf{CZF}_\mathcal{P} \vdash \exists a \in \tf{Ord} \exists u, z \in \tf{L}_a \Big( \varphi_{\ff{nat}}(u) \mathand \forall f \in \prescript{u}{}u \big( f \in \tf{L} \rightarrow f \in z \big) \Big).
\]
By \Cref{theorem:EffectivelyConstructIndex} (\ref{theorem:EffectivelyConstructIndex3}), it follows that we are able to effectively construct an index for an $E_\wp$-recursive function $g$ such that $\tf{CZF}_\mathcal{P} \vdash \forall x \; g(x) \textit{ is defined}$ and
\[
\tf{CZF}_\mathcal{P} \vdash \forall x \; g(x) \rlztp \exists a \in \tf{Ord} \exists u, z \in \tf{L}_a \Big( \varphi_{\ff{nat}}(u) \mathand \forall f \in \prescript{u}{}u \big( f \in \tf{L} \rightarrow f \in z \big) \Big).
\]
Note that $x$ is a dummy variable and the fact that our statement in $\tf{CZF}_\mathcal{P}$ can be realized by the $E_\wp$-recursive function $g$ will allow us to find a bound for the witnessing ordinal $a$. So, setting $x = 0$ we have,
\[
\tf{CZF}_\mathcal{P} \vdash g(0) \rlztp \exists a \in \tf{Ord} \exists u, z \in \tf{L}_a \Big( \varphi_{\ff{nat}}(u) \mathand \forall f \in \prescript{u}{}u \big( f \in \tf{L} \rightarrow f \in z \big) \Big).
\]
Unpacking this by using our definition of realizability with truth, we have that \mbox{$\tf{CZF}_\mathcal{P} \vdash \exists w (w \in g(0))$} and
\[
\tf{CZF}_\mathcal{P} \vdash \forall y \in g(0) \: \pairr y \rlztp \pairl y \in \tf{Ord} \mathand \exists u, z \in \tf{L}_{\pairl y} \Big( \varphi_{\ff{nat}}(u) \mathand \forall f \in \prescript{u}{}u \big( f \in \tf{L} \rightarrow f \in z \big) \Big).
\]
Since realizability implies truth, we have that
\begin{equation} \label{plug}
\tf{CZF}_\mathcal{P} \vdash \forall y \in g(0) \: \Big( \pairl y \in \tf{Ord} \mathand \exists u, z \in \tf{L}_{\pairl y} \Big( \varphi_{\ff{nat}}(u) \mathand \forall f \in \prescript{u}{}u \big( f \in \tf{L} \rightarrow f \in z \big) \Big) \Big). \tag{$*$}
\end{equation}
Now, the important thing to observe is that $g$ was an $E_\wp$ recursive function, which means that the statement ``$ g(0) \textit{ is defined}$'' is definable by a $\Sigma^{\mathcal{P}}$-sentence. Therefore, by \Cref{theorem:CZFPisPi2Conservative}, 
\[
\tf{IKP}(\mathcal{P})\vdash g(0) \textit{ is defined}.
\]

Next, by \Cref{theorem:SigmaReflectstoBH}, we can effectively construct an ordinal $\sigma < \psi_{\Omega}(\varepsilon_{\Omega+1})$ such that the von Neumann hierarchy exists up to level $\sigma$ and $g(0) \in \tf{V}_\sigma$ holds provably in $\tf{IKP}(\mathcal{P})$ and hence in $\tf{CZF}_\mathcal{P}$. 
Thus, for any $y \in g(0)$, $\pairl y \in \tf{V}_\sigma$. Plugging this into (\ref{plug}) and using that $g(0)$ is inhabited, we obtain,
\[
\tf{CZF}_\mathcal{P} \vdash \exists a \in \tf{V}_\sigma \: \Big( a \in \tf{Ord} \mathand \exists u, z \in \tf{L}_a \Big( \varphi_{\ff{nat}}(u) \mathand \forall f \in \prescript{u}{}u \big( f \in \tf{L} \rightarrow f \in z \big) \Big) \Big).
\]
Thus, there is some $a \in \tf{L}_\sigma$ for which $\prescript{\omega}{} \omega \cap \tf{L}$ is an element of $\tf{L}_a$. But from this it follows that $\prescript{\omega}{} \omega \cap \tf{L}_a = \prescript{\omega}{} \omega \cap \tf{L}_\sigma$, so
\[
\tf{CZF}_\mathcal{P} \vdash \exists a \in \tf{L}_\sigma \: \forall f \in \prescript{\omega}{}\omega \big( f \in L_\sigma \leftrightarrow f \in \tf{L}_a \big).
\]
However, this is now a $\Sigma$-sentence, since the $\tf{L}_\alpha$ hierarchy is defined by $\Sigma$-recursion, which means that the same sentence holds in $\tf{IKP}(\mathcal{P})$ by \Cref{theorem:CZFPisPi2Conservative}. Therefore, the result will also hold in the  theory $\tf{KP}(\mathcal{P})$. We say that an ordinal $\alpha$ is a \emph{gap ordinal} if $(\tf{L}_{\alpha+1} \setdiff \tf{L}_\alpha) \cap \mathcal{P}(\omega) = \emptyset$ and clearly the ordinal $a$ above must be a gap ordinal.
 However, as proven in \cite{LeedsPutnamGapOrdinals}, the first gap ordinal is equal to the smallest ordinal such that $\tf{L}_\beta \models \tf{ZFC}^{\uminus}$ (where $\tf{ZFC}^{\uminus}$ stands for $\tf{ZFC}$ without powerset) an ordinal much larger than the Bachmann-Howard ordinal. Thus we have obtained our desired contradiction. 
\end{proof}

We will not explicitly state the background theory needed for the proof of the above result but remark here that it will follow from $\tf{KP}(\mathcal{P})$ augmented by the statement that the von Neumann hierarchy $V_\alpha$ exists for every $\alpha \leq \psi_{\Omega}(\varepsilon_{\Omega+1})$. With more analysis this could be significantly weakened.

\section{Constructive Ordinals}

As mentioned at the beginning of this paper, being an inner model over a classical theory has three defining properties:
\begin{itemize}
\item It is a subclass of the universe;
\item It models the same base theory as the universe;
\item It contains every ordinal in the universe.
\end{itemize}
While we have proven that \tf{L} is a definable subclass of the universe and if the universe satisfies one of $\tf{IKP}^{\minusinf}$, \tf{IKP} or \tf{IZF} then so does \tf{L}, we have also shown that the same cannot be said for \tf{CZF} and $\tf{CZF}_\mathcal{P}$. Moreover, there is one significant, fundamental and structural question that remains open,

\begin{question}[Lubarsky, \cite{LubarskyIntuitionisticL}]
Does \tf{L} contain every ordinal? That is, is $\tf{L} \cap \tf{Ord} = \tf{V} \cap \tf{Ord}$?
\end{question}

\noindent This appears to be a difficult question to resolve. By \Cref{thm: approx in L} we know that \tf{L} always contains a proper class of ordinals and, as a transitive structure, if \tf{L} contains one ordinal, then it contains every smaller ordinal. However, if Excluded Middle were to fail, then the ordinals can no longer be linearly ordered so while we may have no gaps in the ordinals, our ``tree'' of ordinals could be missing branches from the ground model. It is worth remarking that \tf{L} does at least contain every \emph{truth value} in \tf{V}.

\begin{prop}
For every $\alpha \in \Omega = \{ x \divline x \subseteq 1 \}$, $\alpha = \tf{L}_\alpha \in \tf{L}$.
\end{prop}

\begin{proof}
Suppose $\alpha \in \Omega$, then $\tf{L}_\alpha \coloneqq \bigcup_{\beta \in \alpha} \defop(\tf{L}_\beta) = \{ x \divline \exists \beta \in \alpha (x \in \defop(\tf{L}_\beta)) \}$. Now, since $\alpha \subseteq 1$, $\tf{L}_\alpha = \{ x \divline \exists \beta \in 1 (\beta \in \alpha \wedge x \in \defop(\tf{L}_\beta)\}$. But if $\beta \in 1$ then $\defop(\tf{L}_\beta) = \{0\}$ and therefore we must have that $\tf{L}_\alpha = \{ x \divline 0 \in \alpha \wedge x = 0 \} = \{ x \divline x \in \alpha \} = \alpha$.
\end{proof}

\noindent A difficulty arises in trying to prove that \tf{L} contains infinite ordinals which are somehow produced by functions which are not definable in \tf{L}. As a final motivating example for future work, suppose that  we were in an intuitionistic setting with $\tf{V} \neq \tf{L}$ and $c \colon \omega \rightarrow 2$ was some function that was not in \tf{L}. Let $\alpha \in \Omega$ be an arbitrary truth value (which we cannot prove to be either $0$ or $1$) and let
\[
\delta_c \coloneqq \bigcup_{n \in \omega} ((n+2) \cup \{\alpha\}) + c(n).
\]
As a union of ordinals, $\delta_c$ is again an ordinal however it does not seem possible to show that $\delta_c \in \tf{L}$ because to construct it one somehow needs to know about the non-constructible function $c$. Moreover, when one looks at the ordinal $(\delta_c)^\star$ produced in \Cref{thm: approx in L} such that $\tf{L}_{\delta_c} = \tf{L}_{(\delta_c)^\star}$, the construction includes $((n+2) \cup \{\alpha\}) + 1$ for every $n \in \omega$ which erases the function $c$ from the definition. On the other hand, it is difficult to construct an explicit example for which we can prove that we cannot deduce that $\delta_c \in \tf{L}$. The main strategies we have to build intuitionistic models are via realizability, Kripke model or intuitionistic forcing. In each of these cases one adds many new ordinals and we do not know how these strange new ordinals affect the construction of the Constructible Universe.

\begin{question}
Under the above construction, is it always the case that $\delta_c$ is provably in \tf{L}?
\end{question}

\medskip

\subsection*{Acknowledgements} 
\noindent Some of this paper formed part of the first author's PhD thesis, \cite{MatthewsPhD}, which was supervised by the second author and Andrew Brooke-Taylor. The first author was supported by the UK Engineering and Physical Sciences Research Council at the University of Leeds while writing this paper and is grateful for their support during this research. The second author has been supported by the John Templeton Foundation (Grant ID 60842) and is grateful for their support.

We also thank the anonymous referees for their very helpful and insightful comments.

\bibliographystyle{elsarticle-harv}
\bibliography{Constructing_L}

\end{document}